\documentclass{article}

\usepackage{amssymb,amsmath,latexsym,amsthm}

\newfont{\Bb}{msbm10 scaled\magstep0}
\newfont{\Bbs}{msbm10 scaled 800}

\newcommand{\pr}{\mbox{\Bb{P}}}
\newcommand{\fq}{\mbox{\Bb{F}}_q}
\newcommand{\fqs}{\mbox{\Bb{F}}_{q^s}}
\newcommand{\fqbar}{\bar{\mbox{\Bb{F}}}_q}
\newcommand{\fp}{\mbox{\Bb{F}}_p}
\newcommand{\Ql}{\mbox{\Bb{Q}}_\ell}
\newcommand{\Qp}{\mbox{\Bb{Q}}_p}
\newcommand{\Zp}{\mbox{\Bb{Z}}_p}

\newcommand{\Z}{\mbox{\Bb{Z}}}
\newcommand{\sZ}{\mbox{\Bbs{Z}}}
\newcommand{\Q}{\mbox{\Bb{Q}}}

\newcommand{\gX}{{\mathcal X}}
\newcommand{\gS}{{\mathcal S}}
\newcommand{\Oh}{{\mathcal O}}
\newcommand{\E}{{\mathcal E}}
\newcommand{\B}{{\mathcal B}}
\newcommand{\R}{\mbox{\Bb{R}}}
\newcommand{\cR}{{\mathcal R}}
\newcommand{\V}{{\mathcal V}}

\newcommand{\Tr}{{\rm Tr}}
\newcommand{\ord}{{\rm ord}}
\newcommand{\coker}{{\rm coker}}
\newcommand{\Image}{{\rm Im}}
\newcommand{\Spec}{{\rm Spec}}
\newcommand{\Res}{{\rm Res}}
\newcommand{\res}{\tilde{r}} 
\newcommand{\Adj}{{\rm Adj}}
\newcommand{\Fil}{{\rm Fil}}
\newcommand{\reldim}{{\rm rel.dim}}
\newcommand{\Gr}{{\rm Gr}}
\newcommand{\h}{{\mathcal H}}

\newcommand{\Coeff}{{\rm Coeff}}

\newtheorem{theorem}{Theorem}[section]
\newtheorem{lemma}[theorem]{Lemma}
\newtheorem{proposition}[theorem]{Proposition}
\newtheorem{conjecture}[theorem]{Conjecture}

\newtheorem{definition}[theorem]{Definition}
\newtheorem{problem}[theorem]{Problem}

\newenvironment{exafont}{\begin{bf}}{\end{bf}}
\newenvironment{example}{\vspace{0.3cm}\par\noindent\refstepcounter{theorem}\begin{exafont}Example
    \thetheorem\end{exafont}\hspace{\labelsep}}{\vspace{0.3cm}\par}

\newenvironment{note}{\vspace{0.3cm}\par\noindent\refstepcounter{theorem}\begin{exafont}Note
    \thetheorem\end{exafont}\hspace{\labelsep}}{\vspace{0.3cm}\par}

\title{A recursive method for computing\\ zeta functions of varieties}

\author{Alan G.B. Lauder\footnote{Address: Mathematical
Institute, 24-29 St Giles, Oxford OX1 3LB, UK. E-mail: lauder@maths.ox.ac.uk. The author is a Royal Society University Research Fellow}}

\begin{document}

\maketitle

\section{Introduction}

We present a method for calculating the zeta function of a smooth
projective variety over a finite field which proceeds by induction on
the dimension. Specifically, we outline an algorithm which reduces the
problem of calculating a numerical approximation for the action of
Frobenius on the middle-dimensional rigid cohomology of a smooth
variety, to that of performing the same calculation for a smooth
hyperplane section. We present in detail the main new algorithmic
ingredient under some simplifying assumptions, and give full details
of our algorithm for calculating zeta functions for some specific
surfaces; we call it the ``fibration algorithm''. 
We have implemented the fibration algorithm for these surfaces over prime fields using
the Magma programming language, and present some explicit examples
which we have computed. 

To illustrate the main idea behind our approach, we begin by outlining
the proof given by Deligne of the Riemann hypothesis for a smooth
projective variety $X$ over the finite field $\fq$ \cite{D}.
Specifically, the statement that for each $0 \leq i \leq 2\dim(X)$ the
action of the Frobenius endomorphism on the $\ell$-adic \'{e}tale
cohomology space $H_{et}^i(X,\Ql)$ has eigenvalues of
complex absolute value $q^{i/2}$.

Let $X \subset \pr$ be a smooth projective variety of dimension $n+1 >
1$ defined over the finite field $\fq$. Denote by $\check{\pr}$ the
dual projective space whose points $t$ correspond to hyperplanes $H_t$
in $\pr$, and let $D$ be a line in $\check{\pr}$. Let $\tilde{X}
\subset X \times D$ denote the set of points $(x,t)$ such that $x \in
H_t$. Projection on the first and second coordinates yields maps $X
\stackrel{\pi}{\leftarrow} \tilde{X} \stackrel{f}{\rightarrow} D$. The
fibre of $f$ at $t \in D$ is the hyperplane section $X_t = X \cap H_t$
of $X$. For sufficiently general $D$ these maps define a {\it
  Lefschetz pencil} \cite[(5.1)]{D} (one may need to change the
projective embedding first \cite[(5.7)]{D}).

The action of the Frobenius endomorphism on the $\ell$-adic \'{e}tale
cohomology $H_{et}^*(X,\Ql)$ may be studied via this Lefschetz pencil.
In particular, assuming the result holds for smooth curves and arguing
by induction on the dimension $n+1$, one can reduce the proof of the
Riemann hypothesis for $X$ to the case of the Frobenius action on the
middle-dimensional cohomology space $H_{et}^{n+1}(X,\Ql)$. The Leray
spectral sequence for $f$ and further inductive arguments now reduce
the proof of the Riemann hypothesis to the case of $E_{2,et}^{1,n} :=
H_{et}^1 (D,R^nf_* \Ql)$ \cite[(7.1)]{D}.  Specifically, one must
prove that the Frobenius acting on this finite dimensional
$\Ql$-vector space has eigenvalues which have complex absolute value
$q^{(n+1)/2}$. This is the ``core
problem'' and it requires considerable ingenuity.

In this article, we are interested in computing the eigenvalues of
Frobenius, rather than proving that they verify Weil's conjecture.
However, it should be possible to bring to bear upon this
computational problem the above geometric machinery. Specifically, one
expects that the geometric techniques which Deligne used in his
reduction to the core problem can be made algorithmic. However, even
once this is done, one is still faced with a difficult problem viz.
calculation of the Frobenius action on $E_{2,et}^{1,n}$. The present
author has no idea on how this might be achieved. However, the sketch
of Deligne's proof can be presented in the terms of {\it rigid
  cohomology}, rather than $\ell$-adic \'{e}tale cohomology, and this
theory is much more amenable to computation.  In this article we
present an algorithmic solution to the analogous ``core problem'', at
least under certain simplifying assumptions. The principal novelty of
this algorithmic technique is that it proceeds by induction on the
dimension. Specifically, the calculation of a matrix for the action of
Frobenius on the rigid cohomological analogue $E_{2,rig}^{1,n}$
requires as input a matrix for the action of Frobenius on
$H_{rig}^n(X_t)$ for some hyperplane section $X_t$ of $X$. So we can
show that for the purposes of computation, the ``core problem'' of
calculating Frobenius in the middle dimension can be efficiently
reduced to that of a single instance of the problem one dimension
lower down.  In our method the base case of curves is handled using
Kedlaya's algorithm \cite{KK1}.

We note in passing that for smooth projective hypersurfaces (of odd
dimension) Deligne's solution of the ``core problem'' can be applied
in a different manner, viz. rather than fibring the hypersurface $X$ in
a Lefschetz pencil, one can embed it as a fibre in such a pencil
\cite[(5.12)]{D},\cite{NK2}.  Such an approach to calculating zeta
functions was taken by the author in the ``deformation algorithm''
\cite{LSH}.  From a computational point of view, this latter approach
has the disadvantage that the ``total space'' under consideration has
dimension one more than the hypersurface itself. This impacts somewhat
on the complexity of the ``deformation algorithm''. Specifically, the
time/space complexity in terms of the ``middle betti number'' $\dim
H^{n+1}_{rig}(X)$ is rather high. Our new approach, of fibring the
original variety, though more complicated, does appear better from the
point of view of complexity dependence on the ``betti numbers''.

The algorithm presented in this paper uses the main
technique developed for the ``deformation algorithm'', combined
with a ``higher rank'' generalisation of Kedlaya's algorithm.
Although our recursive approach was conceived
as a general purpose algorithm, our implementation and complexity analysis for some
surfaces suggest it
is likely to be of most use for surfaces which can be fibered into
low genus curves. Specifically, for the surfaces we consider in Sections
\ref{Sec-AffineSurface}, \ref{Sec-CompactSurface} and \ref{Sec-CompData}, if one
fixes the genus $g$ of the generic fibre of the fibration, then the asymptotic complexity
of our algorithm is quasi-quartic in the middle betti number, with quasi-cubic
space requirement. In fact, the complexity in this case is comparable to that
in the original algorithm of Kedlaya \cite{KK1}, only in this case we have surfaces
rather than curves (Theorem \ref{Thm-Varyh}). The dependence on the
genus $g$ itself is roughly comparable to that in the ``deformation algorithm'' for curves; see the
end of Section \ref{Sec-Complexity}.

We now outline the contents of the various sections in this paper. In
Section \ref{Sec-Def} we define the zeta function of a variety
and explain the computational problem which pertains to them. In
Section \ref{Sec-RC} we give the main definitions from rigid
cohomology which we shall need, and define the specific computational
problem on which we shall focus (Problem \ref{Prob-Main}), viz.,
calculation of Frobenius on the space $E_{2,rig}^{1,n}$. Neither
Section \ref{Sec-Def} nor \ref{Sec-RC} contains any original
contribution. 

Section \ref{Sec-RedAlgs} considers an ``abstract''
version of the main computational problem, and proves a number of
theorems relevant to its solution (Theorems \ref{Thm-FGEdR} and
\ref{Thm-Growth}).  The main theorem stated in this section (Theorem
\ref{Thm-CT}) is not new; however, Theorem \ref{Thm-FGEdR} and
\ref{Thm-Growth} together yield an algorithmic/effective proof of a
slight weakening of Theorem \ref{Thm-CT}.  This algorithmic/effective
proof is a new contribution.  The material in Section
\ref{Sec-RedAlgs} amounts to a special case of a ``higher rank'' generalisation
of Kedlaya's algorithm.
Section \ref{Sec-DefFrob}
contains a description of the main technique used in the ``deformation
algorithm''. The analysis of the loss of numerical precision during the application of this
technique is the only original contribution in this
section; see Theorem \ref{Thm-DiffSysLoss} and the discussion following it.  
Section \ref{Sec-MainAlg} presents our algorithmic solution
to the main computational problem. Specifically, we assemble together
the algorithmic and theoretical techniques developed in Sections
\ref{Sec-RedAlgs} and \ref{Sec-DefFrob} to address Problem
\ref{Prob-Main}.  

Section \ref{Sec-AffineSurface} presents an explicit family of
surfaces, viz., open subsets of affine surfaces defined by equations of
the form $Z^2 = \bar{Q}(X,Y)$ under some smoothness assumptions. (We
note that these surfaces were previously studied for different reasons
by the author in his expository papers \cite{LFFA, LJNTB}; see
also the Ph.D. work of Hubrechts \cite{HH}.) The
algorithm described in Section \ref{Sec-MainAlg}, together with an
auxiliary algorithm (Section \ref{Sec-ComputeNabla}), Kedlaya's
algorithm for hyperelliptic curves, and some propositions
(\ref{Prop-SSdeg}, \ref{Prop-Eff1}, and \ref{Prop-Eff2}) allow the
efficient computation of numerical approximations to the action of Frobenius
on the middle-dimensional rigid cohomology of these open surfaces;
see Theorem \ref{Thm-TimeSpace} for a complexity estimate.
In Section \ref{Sec-CompactSurface} we consider smooth compactifications
of these open surfaces, and describe how one may efficiently compute
the full zeta functions of these compact surfaces using the main result of
Section \ref{Sec-AffineSurface}; see Section \ref{Sec-Complexity} for
complexity estimates.
 We have implemented this zeta function algorithm for the
case in which the base field is prime using the Magma programming
language. Section \ref{Sec-CompData} presents some explicit zeta
functions we have computed using our implementation.  

The author would like to make a comment regarding the original
motivation of this work: An interesting problem when calculating zeta
functions using rigid cohomology is establishing good bounds on the
loss of numerical precision, i.e., quantifying the divisions by the
characteristic $p$ which occur during the algorithms. It was the
author's attempt to prove such precision-loss bounds by induction on
the dimension using a deep theorem of Christol-Dwork (see \cite{CD} or
\cite[Chap. V]{DGS}) which lead him to consider a recursive approach
to computing zeta functions. The Christol-Dwork theorem, which can be
thought of as a special case of an {\it effective} $p$-adic local
monodromy theorem, remains an essential ingredient in the theoretical
analysis of the algorithm presented in this paper.

{\it Acknowledgements:} The author would like to thank his colleagues
in the Mathematical Institute and Hertford College, Oxford, and his family for their support
and help. He has also been greatly assisted by the kind help of  Francesco Baldassarri, Gilles Christol, Jan Denef, Bas Edixhoven,
Johan de Jong, Ralf Gerkmann,  Kiran
Kedlaya, Michael Singer, Frederik Vercauteren and Daqing Wan. Especial thanks to Nobuo Tsuzuki for his
detailed personal communication, and permission to include part of it as Section \ref{Sec-Tsuzuki}
in this paper.

\section{Varieties and zeta functions}\label{Sec-Def}

Let $\fq$ be the field finite with $q$ elements of characteristic $p$, and fix
an algebraic closure $\fqbar \supset \fq$. For each integer $s \geq 1$, let
$\fqs$ denote the unique subfield of $\fqbar$ of order $q^s$.
Let
$X$ be a variety defined over $\fq$. For $s \geq 1$,
let $|X(\fqs)|$ be the number of $\fqs$-rational points on $X$.

\begin{definition}
The zeta function of $X$ is the formal power series
\[ Z(X,T) := \exp\left(\sum_{s = 1}^\infty \frac{|X(\fqs)|}{s} T^s \right).\]
\end{definition}

\begin{theorem}[Dwork]
The zeta function is a rational function. Specifically,  $Z(X,T) = P(T)/Q(T)$ for some polynomials
$P(T), Q(T) \in 1 + T \Z[T]$ with $\gcd(P,Q) = 1$. 
\end{theorem}

Suppose that the variety $X$ can be specified using $I(X)$ bits of
data (input size) and that zeta function $P(T)/Q(T)$ requires at most
$O(X)$ bits of data (bound on output size). Let $S(X) =
\max\{I(X),O(X)\}$ (problem size).

The central problem in the ``algorithmic theory of zeta functions''
\cite{DW} is :

\begin{problem}\label{Prob-1}
Given $X$ compute $Z(X,T)$ in a polynomial number of bit operations in
$S(X)$.
\end{problem}

The proof of Dwork's rationality theorem can be transformed into an
algorithm for Problem \ref{Prob-1}, see \cite{LW}.  For a hypersurface
the running time in bit operations is polynomial in $(p
S(X))^{\dim(X)}$, i.e., it solves Problem 1 for hypersurfaces assuming
the dimension is fixed and the characteristic ``small''. This
algorithm though is of little practical interest. The algorithm of
Schoof-Pila solves Problem \ref{Prob-1} for smooth plane projective
curves of fixed degree \cite{JP,RS}; this is the most general result obtained so
far using the $l$-adic theory.

Let us for the remainder of this section assume that $X$ is smooth. Then
\[ Z(X,T) = \prod_{i = 0}^{2\dim(X)} \det(1 - T q^{\dim(X)} F_q^{-1}|
H^i_{rig}(X))^{(-1)^{i+1}}.\] The rigid cohomology spaces
$H^i_{rig}(X)$ are finite dimensional vector spaces over $K$, the
unramified extension of $\Qp$ of degree $e := [\fq:\fp]$. The linear
map $F_q$ is that induced by functorality from the $q$th power
Frobenius endomorphism of the coordinate rings of an affine cover of
$X$. We have $F_q = F_p^e$ where $F_p$ is semi-linear with respect to
the Frobenius endomorphism of the $p$-adic field $K$. We denote $F :=
F_p$. The central problem in the ``algorithmic theory'' of rigid
cohomology is:

\begin{problem}\label{Prob-2}
  For each $i$, compute a sufficienty good numerical approximation to
  a matrix for the semi-linear map $F:H^i_{rig}(X) \rightarrow
  H^i_{rig}(X)$. Moreover, do so in a number of bit operations
  polynomial in $p S(X)$.
\end{problem}

By ``sufficiently good'' we mean good enough to recover the integer
polynomials $P(T)$ and $Q(T)$ from their numerical ($p$-adic)
approximations.

In principle, Kedlaya's algorithm \cite{KK1} can be applied to Problem
\ref{Prob-2}, but the running time of this approach is polynomial in
$(pS(X))^{\dim(X)}$. However, it is a remarkably useful algorithm for
the case $\dim(X) = 1$ where it has been extensively studied and
implemented, see \cite{KKsur}.  Problem \ref{Prob-2} was solved for
smooth projective hypersurfaces in \cite{LSH} (the``deformation
algorithm'') using relative rigid cohomology.  This approach again
seems to be of some practical interest, see
in particular the recent work of Gerkmann \cite{G2} and Hubrechts \cite{HH}. (We refer the
reader to Tsuzuki \cite{NT} for a different approach which also uses
relative rigid cohomology. This method, which one might call the
``degeneration algorithm'', is conceptually very nice; however, it has
only been worked out in one special case and it is not clear to the
present author how widely it can be applied.)

To understand better the performance of the ``deformation algorithm''
it is necessary to look more carefully at the dependence on
input/output size. Specifically, one can consider separately the
dependence on the {\it arithmetic size}, measured by $\log_p(q)$ and
$p$, and the {\it geometric size}.  In the case of smooth projective
hypersurfaces, the latter is the middle betti number $h_2$, say, which is
approximately $(d-1)^n$ where $d$ is the
degree of the hypersurface and $n$ is the dimension. The ``deformation
algorithm'' has good dependence on the arithmetic size. However, the time/space
dependence on the geometric size is rather high. Specifically, based on
the analysis in \cite{HH}, the author conjectures that the ``deformation
algorithm'' requirements $\tilde{\Oh}(p\log(q)^3 h_2^{4+ \omega})$ bit operations
and $\tilde{\Oh}(p\log(q)^3 h_2^5)$ bits of space for a smooth hypersurface
with middle betti number $h_2$ defined over $\fq$. Here the Soft-Oh notation
ignores logarithmic factors \cite[Def. 25.8]{vGG}, and $\omega$ is the exponent for matrix multiplication.

 The aim of the new
approach in this paper is to try and reduce the space/time dependence
on the geometric size by using a more economical geometric method.
This is achieved for the surfaces studied in Section \ref{Sec-AffineSurface}
and \ref{Sec-CompactSurface}; see Section \ref{Sec-Complexity}.

We conclude this section by mentioning two very recent advances in the area:
First, work by Kedlaya et al \cite{AKR} on bounding Picard numbers using
$p$-adic cohomology. Second, forthcoming work of Edixhoven and collaborators on computing
coefficients of certain modular forms using $l$-adic cohomology of high dimensional
varieties.

\section{Rigid cohomology}\label{Sec-RC}

Many authors have contributed to the theoretical
development of rigid cohomology, most notably Berthelot, Dwork and Monsky, and more recently
Kedlaya and Tsuzuki.
We follow the definitions given in Gerkmann \cite[Sec. 3]{G2}
and refer the reader to that source for further
details. 

\subsection{Relative rigid cohomology}\label{Sec-RRC}

Let $k = \fq$ be a finite field of characteristic $p$, and $K$ be the
unramified extension of $\Qp$ of degree $[k:\fp]$. Let $\Oh_K$ denote
the valuation ring of $K$. Then $(p)$ is the maximal ideal of $\Oh_K$
and $\Oh_K/(p) \cong \fq$. Let $\ord_p$ denote the $p$-adic valuation
on $K$ normalised so that $\ord_p(p) = 1$, and $|\cdot|_p :=
p^{-\ord_p(\cdot)}$ the corresponding norm. Extend the norm and valuation
to polynomial rings and finite dimensional vector spaces over $K$ in the obvious manner.

Let $X$ be a $k$-scheme of finite type. Let $(X,\bar{X},\hat{\gX})$ be
an $\Oh_K$-triple for $X$, viz. an open immersion $j:X \hookrightarrow
\bar{X}$ into a proper $k$-scheme, and an ``admissible'' embedding
$i:\bar{X} \rightarrow \hat{\gX}$ into a formal $\Oh_K$-scheme. For
$S$ a $k$-scheme and $(S,\bar{S},\hat{\gS})$ an $\Oh_K$-triple for
$S$, a morphism $(X,\bar{X},\hat{\gX}) \rightarrow
(S,\bar{S},\hat{\gS})$ is a commutative diagram
\[
\begin{array}{rcrcl}
X & \hookrightarrow & \bar{X} & \rightarrow & \hat{\gX}\\
f \downarrow &  & \bar{f}\downarrow & & \downarrow \hat{f}\\
S & \hookrightarrow & \bar{S} & \rightarrow & \hat{\gS}.
\end{array}
\]
The {\it relative rigid cohomology} sheaf of the morphism $f: X \rightarrow S$ is
\[ \h^i_{rig}(X/S) := R^i \hat{f}_{K*} j^\dagger
\Omega^\bullet_{]\bar{X}[/]\bar{S}[}.\]
We take global sections to give the {\it relative rigid cohomology} spaces $H^i_{rig}(X/S) := \Gamma(] \bar{S} [,\h^i_{rig} (X/S))$ with which we shall work.

\subsection{A relative comparison theorem}\label{Sec-RCT}

Let $X$ and $S$ be $k$-schemes of finite type.
Assume now that there exist commutative diagrams
\[
\begin{array}{ccccccc}
  X & \hookrightarrow & \bar{X} &  & S & \hookrightarrow & \bar{S}\\
  \uparrow & & \uparrow & \mbox{and}& \uparrow & & \uparrow\\
  \gX & \hookrightarrow & \bar{\gX} & & \gS & \hookrightarrow & \bar{\gS}.
  \end{array}
\]
Here $\gX,\bar{\gX},\gS,\bar{\gS}$ are $\Oh_K$-schemes. The vertical
maps take the special fibres. Assume that the lower horizontal maps
are open immersions and their codomains $\bar{\gX}$ and $\bar{\gS}$
are proper and smooth $\Oh_K$-schemes.  Now take $\hat{\gX}$ and
$\hat{\gS}$ to be the $p$-adic completions of $\bar{\gX}$ and
$\bar{\gS}$, respectively. Then we obtain $\Oh_K$-triples
$(X,\bar{X},\hat{\gX})$ and $(S,\bar{S},\hat{\gS})$. Assume now that
we have morphisms $f:X \rightarrow S$ and $\bar{f}:\bar{X} \rightarrow
\bar{S}$ as in Section \ref{Sec-RRC}, and further morphisms $\gX
\rightarrow \gS$ and $\bar{\gX} \rightarrow \bar{\gS}$. Assume
moreover that these fit together with the two diagrams immediately
above to give a ``commutative cube''. (By ``completion'' this
yields a morphism of $\Oh_K$-triples as before.) One may define from
this commutative cube the {\it relative de Rham cohomology} sheaf of the
induced morphism on the generic fibres $f_K:\gX_K \rightarrow \gS_K$:
\[ \h_{dR}^i (\gX_K/\gS_K) := R^i f_{K*} \Omega^\bullet_{\gX_K/\gS_K}.\]
 Taking global sections we define the {\it relative de Rham cohomology} spaces \linebreak
  $H_{dR}^i (\gX_K/\gS_K) := \Gamma(\gS_K, \h_{dR}^i
 (\gX_K/\gS_K))$ with which we shall work.
 
Finally, assume that the complement $\bar{\gX} - \gX$ has smooth
components with normal crossings over $\bar{\gS}$. 
Then $\h^i_{dR}(\gX_K/\gS_K)$ is coherent.\footnote{This result
appears to be ``well-known to experts'', although the author has not been able to
find an explicit reference for it. The point is that in such a case $\h^i_{dR}(\gX_K/\gS_K)$ may be computed
via the hypercohomology of a proper morphism applied to a relative logarithmic de Rham complex,
which is coherent. Note that in our application in Section \ref{Sec-AffineSurface} we shall prove finiteness
directly.}
Moreover, according to
Gerkmann \cite[Eqn (8)]{G2}, the comparison theorem of
Baldassarri-Chiarelletto extends to this relative situation.
Specifically, the natural morphism
\begin{equation}\label{Eqn-BCG}
 \h_{dR}^i (\gX_K/\gS_K) \otimes_{\Oh_{\gS_K}} j^\dagger
\Oh_{]\bar{S}[} \rightarrow \h_{rig}^i(X/S)
\end{equation}
is an isomorphism.
Define $A^\dagger := \Gamma(]\bar{S}[, j^\dagger \Oh_{]\bar{S}[})$ and
$A := \Gamma(\gS_K,\Oh_{\gS_K})$ to be the rings
of global sections. Then \cite[Prop. (2.5.2)(ii)]{PB} shows that (\ref{Eqn-BCG})
implies the following theorem.

\begin{theorem}\label{Thm-BCG}
With the assumptions stated in this section, the following isomorphism holds:
\[ H_{dR}^i (\gX_K/\gS_K) \otimes_{A} A^\dagger \cong H^i_{rig}(X/S).\]
\end{theorem}

\subsection{Proper and smooth base change}\label{Sec-BC}

We retain the definitions and assumptions from Section \ref{Sec-RRC} and define
$A^\dagger := \Gamma(]\bar{S}[, j^\dagger \Oh_{]\bar{S}[})$.
Assume now that the morphism $\hat{\gX} \rightarrow \hat{\gS}$ is
proper and smooth. For each point $\bar{\gamma} \in S$ in the base
denote by $X_{\bar{\gamma}}$ the fibre $X \rightarrow S$ at
$\bar{\gamma}$.  The following base change theorem will be of
importance to us \cite[Thm 3.1]{G2}:

\begin{theorem}[Berthelot] \label{Thm-Bert}
Base change $\Spec(k(\bar{\gamma})) \rightarrow S$ induces an isomorphism
$H_{rig}(X/S) \otimes_{A^\dagger} K(\gamma) \cong H_{rig}(X_{\bar{\gamma}})$.
\end{theorem}

\subsection{Pencils of varieties}\label{Sec-Pencils}

We retain the definitions and assumptions in Section \ref{Sec-RCT} and
\ref{Sec-BC}, i.e., we have a morphism $f:X \rightarrow S$, along with
all the auxiliary objects and properties so that the comparison
(Theorem \ref{Thm-BCG}) and base change (Theorem \ref{Thm-Bert})
theorems hold.

Assume now that 
\[ S  = \pr_k^ 1 - \{\bar{\gamma}_1,\dots,\bar{\gamma}_d,\infty\} \]
where the $\bar{\gamma}_i \in k$ are distinct.  Thus $\bar{S} = \pr_k^1$. Choose
$\gamma_i \in \Oh_K$ so that
$\gamma_i \bmod{p} = \bar{\gamma}_i$. Take
\[ \gS = \pr_{\Oh_K}^ 1 - \{\gamma_1,\dots,\gamma_d,\infty\}. \] The
coordinate ring $A := \Gamma(\gS_K,\Oh_{\gS_K})$ of the generic fibre $\gS_K$ is the
localisation $K[\Gamma][1/r(\Gamma)]$ where $r := \prod_{i = 1}^d
(\Gamma - \gamma_i) \in K[\Gamma]$.  The ring of global sections
$A^\dagger := \Gamma(]\bar{S}[,j^\dagger\Oh_{]\bar{S}[})$ is the weak (a.k.a. dagger) completion of
$K[\Gamma][1/r(\Gamma)]$. We denote this ring by
$K[\Gamma][1/r(\Gamma)]^\dagger$. Its elements can be written in the
form $\sum_{i \in \sZ} a_i(\Gamma)r(\Gamma)^i$ where the coefficients
$a_i(\Gamma) \in K[\Gamma]$ have $\deg(a_i) < \deg(r)$ and satisfy
$\ord_p(a_i) - \epsilon |i| \rightarrow \infty$ as $|i| \rightarrow
\infty$ for some $\epsilon > 0$. (Note that there is no lower bound on
the $\epsilon$ which occur for different elements in the ring.)  It is
convenient at this stage to give a definition which we shall need
later.

\begin{definition}\label{Def-EPB}
  We shall say that we have given {\rm effective $p$-adic bounds} for
  an element in $a = \sum_{i \in \sZ} a_i(\Gamma) r(\Gamma)^i \in
  K[\Gamma,1/r(\Gamma)]^\dagger$ if we are given $\eta,\delta \in \Q$
  with $\eta > 0$ such that $\ord_p(a_i) \geq \eta |i| + \delta$ for
  all $i \in \Z$.
\end{definition}

\subsection{Cohomology in the middle dimension}\label{Sec-CMD}

We retain the definitions and assumptions from Section
\ref{Sec-Pencils}, i.e., we have a pencil $f:X \rightarrow S$ and the
comparison and base change theorems hold.

Let us now focus our attention on the middle dimension. Specifically,
let $n$ be the relative dimension of the morphism $f:X \rightarrow S$.
By the comparison theorem and coherence of relative de Rham
cohomology, 
 we see that
$H^n_{dR}(\gX_K/\gS_K)$ and $H^n_{rig}(X/S)$ are locally free modules
of finite rank over the rings $A = K[\Gamma,1/r(\Gamma)]$ and
$A^\dagger = K[\Gamma,1/r(\Gamma)]^\dagger$, respectively.  We shall {\it assume}
that they are in fact free. Let us simplify our notation now by writing
\[ \E := H^n_{dR}(\gX_K/\gS_K),\, \E^\dagger := H^n_{rig}(X/S),\]
so by the comparison theorem (Theorem \ref{Thm-BCG}) we have 
\[ \E^\dagger = \E \otimes_A A^\dagger.\] The free modules $\E$ and
$\E^\dagger$ come with additional structure. Specifically, derivation
with respect to $\Gamma$ induces a {\it connection} $\nabla$ on the
$A$-module $\E$. Let us recall the precise definition of a connection.

\begin{definition}\label{Def-Connection}
A connection $\nabla $ on $\E$ is a  map $\nabla:\E\rightarrow \E \otimes \Omega^1_A$ such that
$\nabla(e_ 1 + e_2) = \nabla(e_1) + \nabla(e_2)$ and $\nabla(ae_1) = 
e_1 \otimes da + a \nabla(e_1)$ for all  $a \in A$ and $e_1,e_2 \in \E$.
\end{definition}

Here $d:A \rightarrow \Omega^1_{A}$ is the universal derivation, which
in our case amounts to differentiation w.r.t. $\Gamma$.  The
connection induced by differentiation w.r.t. $\Gamma$ is called the
{\it Gauss-Manin} connection, and the pair $(\E,\nabla)$ a {\it
  $\nabla$-module}. Differentiation with respect to $\Gamma$ also
induces a connection $\nabla^\dagger:\E^\dagger \rightarrow \E^\dagger
\otimes \Omega^1_{A^\dagger}$; here $\Omega^1_{A^\dagger}$ is the
module of $K$-linear differentials which are continuous w.r.t. the
$p$-adic norm.  The comparison theorem tells us that this is just the
Gauss-Manin connection. Specifically, assuming one can compute a
matrix for $\nabla$ with respect to some basis of $\E$ over $A$, then
the {\it same} matrix defines the map $\nabla^\dagger$ on the basis of
$\E^\dagger$ obtained by extending scalars. This means that the
Gauss-Manin connection $\nabla^\dagger$ on $\E^\dagger$ can be
computed in a purely algebraic manner.

The module $\E^\dagger$ comes with one further piece of data, namely
the {\it Frobenius} map $F:\E^\dagger \rightarrow \E^\dagger$. This
map is induced by the functorality of the construction from the $p$th
power map on the coordinate rings of an affine cover of $X$. It is
additive and semi-linear. Specifically, the construction requires the
choice of a lifting of the $p$th power map from $A\otimes_K \fq$ to
$A^\dagger$. Let us {\it assume} that we have chosen the obvious
lifting $\sigma:A^\dagger \rightarrow A^\dagger$ so that
$\sigma:\Gamma \rightarrow \Gamma^p$ and $\sigma$ acts on $K$ as the
Frobenius automorphism. Then $F$ is a $\sigma$-linear map on the
module $\E^\dagger$. Thus with respect to a basis for $\E^\dagger$, it
can also be described by a matrix. Note though that this matrix has
entries in the ring $A^\dagger$, whereas the matrix for
$\nabla^\dagger$ referred to in the previous paragraph has entries in
$A$.  Since $\frac{d}{d\Gamma}$ and $\sigma$ commute up to a factor
$p\Gamma^{p-1}$ on the weak completions of the coordinate rings of an affine cover of $\gX_K$,
functorality of the construction yields the following important
commutative diagram:
\begin{equation}\label{Diagram-F}
\begin{array}{rcl}
  \E^\dagger & \stackrel{\nabla^\dagger}{\rightarrow} & \E^\dagger \otimes_A^\dagger \Omega^1_{A^\dagger}\\ 
  F \downarrow & & \downarrow F \otimes d\sigma\\
  \E^\dagger & \stackrel{\nabla^\dagger}{\rightarrow} & \E^\dagger \otimes_A^\dagger
  \Omega^1_{A^\dagger}.
\end{array}
\end{equation}
The data $(\E^\dagger, \nabla^\dagger,F)$
is called a {\it $(\sigma,\nabla^\dagger)$-module} over $A^\dagger$, or alternatively an
{\it overconvergent $F$-isocrystal} over $S$. Let us recall the precise definition.

\begin{definition}\label{Def-Isocrystal}
  An overconvergent $F$-isocrystal $(\E^\dagger,\nabla^\dagger,F)$ on
  $S$ (a.k.a. $(\sigma,\nabla^\dagger)$-module over $A^\dagger$)
  consists of the following data. A (locally) free $A^\dagger$-module
  $\E^\dagger$, an injective $\sigma$-linear map $F:\E^\dagger
  \rightarrow \E^\dagger$, and a connection $\nabla^\dagger:\E^\dagger
  \rightarrow \E^\dagger \otimes \Omega^1_{A^\dagger}$, such that
  Diagram (\ref{Diagram-F}) commutes.
\end{definition}

We shall denote the kernel and cokernel of the map $\nabla^\dagger$ in
Diagram (\ref{Diagram-F}) by $H^0_{rig}(S,\E^\dagger)$ and
$H^1_{rig}(S,\E^\dagger)$, respectively. These objects are vector
spaces over $K$. By commutativity, $F$ induces a map on each of these
spaces. We note that the space
$\ker(\nabla^\dagger) = H^0_{rig}(S,\E^\dagger)$ is certainly
finite dimensional  over $K$. (It embeds in the finite-dimensional space of local solutions around any
non-singular point; an observation of Monsky.)

We now state our final assumption on the family $X \rightarrow S$. We
assume that the family $\gX_K \rightarrow \gS_K$ comes by extension
of scalars from a smooth morphism defined over an algebraic number field. It follows
then by the ``open local monodromy theorem''  that the connection is regular, i.e., locally has only simple
poles, and the local exponents (see Section \ref{Sec-DefEetc}) are
rational numbers \cite[Thm. (14.3)]{NK}.

Under this final assumption, as well as the others already in place in
this section, it is known that $H^1_{rig}(S,\E^\dagger)$ is also
finite dimensional (follows from \cite[Corollary 2]{BC}).  We will
give an effective/algorithmic proof of finiteness under some
simplifying assumptions (follows from Theorems \ref{Thm-FGEdR} and
\ref{Thm-Growth}). (There is also an older, related result due to
Adolphson \cite[Theorem 2, Remark p. 286]{AA}.)

We now come to the main definition in the paper.

\begin{definition}\label{Def-E2}
The space $E_{2,rig}^{1,n} := H^1_{rig}(S,\E^\dagger) = \coker(\nabla^\dagger:\E^\dagger
\rightarrow \E^\dagger \otimes_A \Omega^1_{A^\dagger})$.
\end{definition}

By a result communicated to us by Professor Nobuo Tsuzuki, when $X$ is affine this space is a term
in a spectral sequence for the morphism $X \rightarrow S$. In fact,  for $X$ affine we
have the isomorphism $H^{n+1}_{rig}(X) \cong E_{2,rig}^{1,n}$, see
Eqn (\ref{Eqn-Tsuzuki}) in Section \ref{Sec-Tsuzuki}.

Finally we are able to state the computational problem we consider:

\begin{problem}\label{Prob-Main}
Calculate a numerical approximation to a matrix for the map $F: E_{2,rig}^{1,n} 
\rightarrow E_{2,rig}^{1,n}$.
\end{problem}

We solve this problem under the assumption that we are given as input suitable
numerical approximations to:
\begin{itemize}
\item{A matrix for the connection $\nabla$, {\it and this matrix
has only simple poles (even ``modulo $p$'') and prepared local
exponents}.}
\item{A specialisation of the matrix
for $F:\E^\dagger \rightarrow \E^\dagger$ at a Teichm\"{u}ller point.}
\end{itemize}
We further assume that:
\begin{itemize}
\item{We are given as input effective $p$-adic estimates for the matrix for $F:
\E^\dagger \rightarrow \E^\dagger$ (see Definition \ref{Def-EPB}).}
\end{itemize}
Regarding the first input, ones expects to be able to compute this
matrix efficiently in any concrete application of the method, e.g., in
the case in which a basis of forms for $\E$ is known it can be done
using linear algebra. See for example our calculation in Section
\ref{Sec-ComputeNabla}. The assumption on the matrix for the
connection is ``locally'' true by the regularity and local monodromy
theorem \cite[Thm. (14.3)]{NK}, since the family $\gX_K \rightarrow \gS_K$ can be
defined over an algebraic number field.  Our {\it simplifying
  assumption} is that there is a {\it global basis} for which the
matrix has only simple poles with prepared local exponents.

Regarding the second input, by the base change theorem (Theorem
\ref{Thm-Bert}) the specialisation, say at $\Gamma = \gamma$ a
Teichm\"{u}ller point, is precisely the matrix for the $p$th power
Frobenius map acting on the cohomology space
$H^n_{rig}(X_{\bar{\gamma}}/\fq(\bar{\gamma}))$.  Here
$X_{\bar{\gamma}}$ is the fibre of the family at $\bar{\gamma} :=
\gamma \bmod{p}$. Such a matrix can be computed recursively (or by
Kedlaya's algorithm in the case $n = 1$).

Regarding the third assumption, again one expects to be able to calculate
such bounds in any concrete application of the method, see Section \ref{Sec-EffBounds}.

\subsection{Leray Spectral Sequence}\label{Sec-Tsuzuki}

This section is independent of the rest of the paper. We describe the contents of a personal
communication from Professor Nobuo Tsuzuki to the author. Note that the notation
in this section is consistent with \cite{NTss}, but varies slightly from that in the
remainder of this paper. 

Let $K$ be a complete discrete valuation field of mixed characteristic $(0,p)$ and
$\V$ and $k$ be the ring of integers and residue field of $K$, respectively. Let
$S$ be an affine smooth scheme of dimension $m$ over $\Spec(k)$ and $X/S$
a smooth family with $n := \reldim(X/S)$ such that $X$ is affine. Suppose there exists
a smooth afffine lift $\gX/\gS$ of $X/S$ over $\Spec(\V)$ with $A = \Gamma(\gX,\Oh_{\gX})$ and
$R = \Gamma(\gS,\Oh_{\gS})$ such that $\Omega^1_{R/\V}$ is a free $R$-module. Then
one can calculate the rigid cohomology of $X/K$ as
\[ H^r_{rig}(X/K) := H^r(A^\dagger_K \otimes \Omega^\bullet_{A/\V})\,(\, =: H^r_{MW}(X/K)),\]
where $A^\dagger$ is the weak (a.k.a. dagger) completion of $A$ over $\V$ and $A^\dagger_K
 := A^\dagger \otimes_\V K$. Let us define a filtration $\Fil^*$ of $A^\dagger_K \otimes_A
 \Omega^\bullet_{A/\V}$ by
 \[ \Fil^q := \Image(A^\dagger_K \otimes_A \Omega^{\bullet - q}_{A/\V} \otimes_R
 \Omega^q_{R/\V} \rightarrow A^\dagger_K \otimes_A \Omega^\bullet_{A/\V}).\]
 Since $\Omega^q_{R/\V}$ is a free $R$-module, one has
 \[ \Gr_{\Fil}^q = A^\dagger_K \otimes_A \Omega^{\bullet - q}_{A/R} \otimes_R \Omega^q_{R/\V}.\]
 There exists a spectral sequence \cite[Thm. 3.4.1]{NTss}
 \[ E_1^{q,r} := H^r(A^\dagger_K \otimes_A \Omega^\bullet_{A/R}) \Rightarrow
 H^{q+r}(A^\dagger_K \otimes_A \Omega^\bullet_{A/\V}) = H^{q+r}_{MW}(X/K),\]
where the edge homomorphism is called the Gauss-Manin connection. Since
$E_1^{q,r} = 0$ except when $0 \leq q \leq m$ and $0 \leq r \leq n$, one has
\[ E_2^{m,n} = \dots = E_\infty^{m,n} = H^{m+n}_{MW}(X/K).\]
Hence the top rigid cohomology group $H_{rig}^{m+n}(X/K)$ is calculated
by the Gauss-Manin connection:
\begin{equation}\label{Eqn-Tsuzuki}
 H_{rig}^{m+n}(X/K) \cong \coker( H^n(A^\dagger_K \otimes_A \Omega^\bullet_{A/R})
\otimes_R \Omega^{m-1}_{R/\V} \rightarrow H^n(A^\dagger_K \otimes_A \Omega^\bullet_{A/R})
 \otimes_R \Omega^m_{R/\V}).
 \end{equation}

\section{Algorithms for reduction in  $E_{2,rig}^{1,n}$}\label{Sec-RedAlgs}

This section is independent of Section \ref{Sec-RC}, although relies
on it for motivation.  We recall the definitions we shall need in an
abstract manner, stripped of their geometric origin.

\subsection{Definitions}\label{Sec-DefEetc}

Let $k=\fq$ be the finite field with $q$ elements of characteristic
$p$, and $K$ the unramified extension of $\Qp$ of degree $[k:\fp]$.
Denote by $\bar{K}$ an algebraic closure of $K$.
Let $\Oh_K$ be the ring of integers of $K$, and $r(\Gamma) \in
\Oh_K[\Gamma]$ a monic polynomial of degree $d$ which is squarefree
modulo $p$. Let $A := K[\Gamma,1/r(\Gamma)]$ and let $A^\dagger$ be
the dagger completion of $A$ (this is defined in Section
\ref{Sec-Pencils}). Let $\E$ be a free module of finite rank $m$ over
$A$ and define $\E^\dagger := \E \otimes_A A^\dagger$. Let $\nabla:\E
\rightarrow \E \otimes \Omega^1_A$ be a connection (Definition
\ref{Def-Connection}).  Fix a basis $\B$
for $\E$ over $A$ and represent elements in $\E$ as column vectors
w.r.t. this basis. Take for $\Omega^1_A$ the basis element $d\Gamma$
over $A$.  Take the basis for $\E \otimes \Omega^1_A$ to be the tensor
product of these two bases.  Assume that with respect to this choice,
the connection $\nabla$ acts as
\begin{equation}\label{Eqn-DefNabla}
  \nabla = \frac{d}{d\Gamma} + \frac{b(\Gamma)}{r(\Gamma)}: \E \cong A^m \rightarrow \E \otimes
  \Omega_A^1 \cong A^m \otimes d\Gamma
 \end{equation}
where the matrix $b(\Gamma) \in M_m(\Oh_K[\Gamma])$ has degree in $\Gamma$ at most
$d - 1$. This assumption ensures that the matrix for $\nabla$ has only simple
poles, including at infinity. This is our {\it main} simplifying assumption. Such a
differential system is called {\it fuchsian}. Any differential systems may in principle, after a change of basis, be written in this form, possibly at the expense of introducing one new pole.
See the discussion of the Riemann-Hilbert problem in \cite[Section 5.3]{vPS}.

Let $\nabla^\dagger:\E^\dagger \rightarrow \E^\dagger \otimes
\Omega^1_{A^\dagger}$ be obtained from $\nabla$ by extension of
scalars. We are interested in the spaces
\[ E_{2,dR}^{1,n} :=\coker(\nabla:\E \rightarrow \E \otimes \Omega^1_A),\, E_{2,rig}^{1,n} :=
\coker(\nabla^\dagger:\E^\dagger \rightarrow \E^\dagger \otimes \Omega^1_{A^\dagger}).\]
The notation chosen here is to remind the reader of the ``geometric origin'' of the
connections we shall actually be considering.

For $R := \{\gamma \in \Oh_{\bar{K}}\,|\,r(\gamma) = 0\}$ note that
\[ \frac{b(\Gamma)}{r(\Gamma)} = \frac{b(\Gamma)}{r^\prime(\Gamma)}
\sum_{\gamma \in R} \frac{1}{\Gamma - \gamma}\] 
where ``dash''
indicates differentiation w.r.t. $\Gamma$. Thus the {\it residue
  matrix} at the regular singular point $\Gamma = \gamma \in R$ is
$b(\gamma)/r^\prime(\gamma)$; the set of eigenvalues of this matrix,
denoted $E_\gamma$, are the {\it local exponents} or {\it local monodromy eigenvalues}
at $\Gamma =
\gamma$. One checks that the residue matrix at infinity is $-b_{d-1}$,
the negative of the coefficient of $\Gamma^{d-1}$ in $b(\Gamma)$.  The set $E_\infty$
is defined as the set of eigenvalues of $-b_{d-1}$. Finally, the {\it
  exponent} set of $\nabla$ w.r.t. the basis $\B$ is
\[ E(\nabla,\B) := E_\infty \cup \bigcup_{\gamma \in R} E_\gamma.\]
Note that this set  modulo $\Z$ is independent of the basis $\B$.

\begin{definition}\label{Def-Rho}
Let $\rho = \rho(\nabla,\B)$ be the smallest positive integer larger than any
integer in the set $E(\nabla,\B)$.
\end{definition}

Denote by $A \otimes d\Gamma_\rho$ the $K$-vector space of $1$-forms
spanned by the set
\[ \left\{ \frac{\Gamma^i}{r^j} \otimes d\Gamma\,:\, 0 \leq i < d, 1
  \leq j \leq \rho\right\} \cup \left\{ \Gamma^j \otimes d\Gamma\,:\,
  0 \leq j \leq \rho - 2\right\}.\] 
Denote by $\E \otimes
d\Gamma_\rho$ the $K$-vector space spanned by column vectors in $A^m
\otimes d\Gamma$ whose entries belong to the space $A \otimes
d\Gamma_\rho$.

\subsection{Effective finiteness of $E^{1,n}_{2,dR}$.}

We can now state our first finiteness theorem.

\begin{theorem}\label{Thm-FGEdR}
Let the pair $(\E,\nabla)$ be as defined in Section \ref{Sec-DefEetc}, and
$\rho = \rho(\nabla,\B)$ the positive integer from Definition \ref{Def-Rho} which
depends upon both $\nabla$ and the basis $\B$ for $\E$.
Then $\coker(\nabla)$ is generated over $K$ by the image of the space $\E \otimes d\Gamma_\rho$.
\end{theorem}

\begin{proof}
We shall give an algorithm for writing 
an element $u \in \E \otimes d\Gamma$ in the form $u = \nabla(v) + w$
with $v \in \E$ and $w \in \E \otimes d\Gamma_\rho$.
It proceeds in two stages: First, simultaneous reduction
of the pole orders of $1$-forms at the roots of $r$; Second, reduction of pole
orders at infinity.

Let $U(\Gamma) \in K[\Gamma]^m$, viewed as a column vector. We shall show that
for $\ell \geq \rho$ we have
\begin{equation}\label{Eqn-Redr}
\frac{U}{r^{\ell+1}} \otimes d\Gamma = \nabla\left(\frac{V}{r^\ell}\right) + \frac{W}{r^\ell} \otimes d\Gamma
\end{equation}
for some $V,W \in K[\Gamma]^m$ with $\deg(V) < d$ and 
\[ \deg(W) \leq \max\{\max\{2d-2, \deg(U)\}- d,0\}.\]
Moreover, we shall give a method for computing $V$ and $W$.

We claim that there exists
a unique $V(\Gamma) \in K[\Gamma]^m$ with $\deg(V) < d = \deg(r)$ such that
\[ (-\ell r^\prime I_m + b)V \equiv U \mod{r}.\]
Let us assume this claim is true. Define 
\[ X:= \frac{(-\ell r^\prime I_m + b)V - U}{r} \in K[\Gamma]^m.\]
Then 
\[ \deg(X) \leq \max\{\max\{2d-2, \deg(U)\} - d, 0\}.\]
Define $W := - X - V^\prime$. Then
$\deg(W)$ is bounded as claimed above  
and one checks by direct computation that (\ref{Eqn-Redr}) holds.

It remains to establish the uniqueness, existence and computability of $V$. For this,
we must show that the determinant of the matrix $(-\ell r^\prime I_m + b)$ is a unit
modulo $r$. Now
\[ (-\ell r^\prime I_m + b) = -r^\prime f(\ell,\Gamma)\]
where
\[ f(t,\Gamma) := t I_m - \frac{b(\Gamma)}{r^\prime
(\Gamma)}.\]
Now $r^\prime$ is invertible modulo $r$ since the latter is squarefree. We need
to show $\det(f(\ell,\Gamma))$ is a unit modulo $r = \prod_{\gamma \in R} (\Gamma - \gamma)$, 
i.e., we must show that $\det(f(\ell,\gamma)) \ne 0$ for all 
$\gamma \in R$. But $\det(f(\ell,\gamma)) = 0$ if and only if $\ell$ is an eigenvalue of
the matrix $b(\gamma)/r^\prime(\gamma)$. Since $\ell \geq \rho$ and $\rho$ is larger
than any integer element in the set $\cup_{\gamma \in R} E_\gamma$ the result follows. 

Let $U(\Gamma) \in K[\Gamma]^m$ as before. We shall show that if
\[ \deg(U) - \rho d =: \ell - 1 > \rho - 2\] 
then  
\begin{equation}\label{Eqn-InfPoles}
\frac{U}{r^\rho} 
\otimes d\Gamma = \nabla (V \Gamma^{\ell}) + \frac{W}{r^\rho} \otimes d\Gamma
\end{equation}
for some $V \in K^m$ and $W \in K[\Gamma]^m$ with $\deg(W) \leq \deg(U) - 1$.
 
We shall take local expansions of rational functions around the origin. Put
\[ U = u_{\ell-1} \Gamma^{\ell - 1} + u_{\ell - 2} \Gamma^{\ell - 2} + \dots,\,
\nabla = \frac{d}{d\Gamma} + \left( b_{d-1} \Gamma^{-1} + \dots \right).\]
Here $b_{d-1} \in M_m(K)$ is the coefficient of
the monomial $\Gamma^{d-1}$ in $b(\Gamma)$. 
Let $V \in K^m$ be the element such that
\[ \left(\ell I_m + b_{d-1}\right)V = u_{\ell - 1}.\]
We note that $V$ exists and is unique
by the assumption that the integers in the eigenvalue set $E_\infty$ of $-b_{d-1}$ are all
less than $\rho$, and $\ell \geq \rho$.
By direct computation one checks that
\[ \frac{U}{r^\rho} \otimes d\Gamma - \nabla( V \Gamma^{\ell}) = \frac{W}{r^\rho} \otimes d\Gamma\]
where $\deg(W) \leq \deg(U) - 1$. This concludes the description of the algorithm.

We note that in implementations one
should represent the numerator $U$ in an $r(\Gamma)$-adic expansion. With such a representation,
in the second stage it is more efficient to compute
\[ \frac{A}{r^\rho} \otimes d\Gamma - \nabla\left(V 
\Gamma^{\ell - d \lfloor \ell/d \rfloor}
r^{\lfloor \ell/d \rfloor} \right) \]
with $V \in K^m$ as in the preceding paragraph.
\end{proof}

\begin{theorem}\label{Thm-BasisdR}
The space $E^{1,n}_{2,dR}$ is a finite dimensional $K$-vector space.
\end{theorem}

\begin{proof}
A basis for this space can be computed using linear algebra and Theorem \ref{Thm-FGEdR}.
Specifically, one computes a basis for the cokernel of the $K$-linear map $\nabla:A^m_{\rho-1}
\rightarrow A^m \otimes d\Gamma_\rho$. Here $A^m_{\rho - 1} \subset A^m$ is the $K$-space of
column vectors whose entries have poles of order at most $\rho - 1$. 
\end{proof}

We note that the  dimension of $E^{1,n}_{2,dR}$ can be calculated explicitly.

\subsection{Effective finiteness of $E^{1,n}_{2,rig}$}

In this section we give an effective/algorithmic proof of the finiteness of
$E^{1,n}_{2,rig}$ under certain conditions. First, consider the conditions:\\
\\
{\it (Rat.) The exponent set $E(\nabla,\B)$ contains only rational numbers.}\\
\\
{\it (Prep.Rat.) The exponent set $E(\nabla,\B)$ contains only rational numbers; moreover,
for each $E \in \{E_\gamma\}_{\gamma \in R} \cup \{E_\infty\}$, if
$\lambda_1,\lambda_2 \in E$ then $\lambda_1 - \lambda_2$ is not a positive integer.}\\
\\
If (Prep.Rat.) holds then we say that the local exponents are {\it prepared}.
Certainly $(Prep.Rat.) \Rightarrow (Rat.)$.  Also,
rationality of the exponents for a regular connection is independent of
the basis chosen.
 
 Next, consider the following definition:

\begin{definition}
The connection $\nabla$ is {\rm overconvergent} if it has a basis of local
solutions which converge on the $p$-adic unit disk around the generic point
$t$ with $|t|_p = 1$, c.f. \cite[Chap. III, Sec. 5]{DGS}.
\end{definition}

We now give a alternative characterisation of overconvergence.

\begin{theorem}\label{Thm-Trans}
The connection $\nabla$ is overconvergent if and only if there exists an element
$\gamma \in \Oh_{\bar{K}}$ with $\ord_p r(\gamma) = 0$ such that the differential
system $\frac{d}{d\Gamma} + b(\Gamma)/r(\Gamma)$ has a basis of local solutions on
the $p$-adic open unit disk around $\Gamma = \gamma$. 
\end{theorem}

\begin{proof}
That this is necessary follows by specialising the generic solution matrix at $\Gamma = \gamma$.
For sufficiency, we observe that \cite[Chap IV Prop 5.1]{DGS} allows one to transfer
the convergence on the open unit disk around the point $\Gamma = \gamma$ to 
the same disk around the generic point. (Specifically, change variables so that
$\gamma = 0$, and take ``$\alpha$'' to be the generic point $t$ with $|t|_p = 1$.)
\end{proof}

\begin{definition}
The ``dual'' connection $\check{\nabla}$ is defined to act as
\[ \check{\nabla}: \frac{d}{d\Gamma} - \frac{b(\Gamma)}{r(\Gamma)}:
A^m \rightarrow A^m \otimes d\Gamma.\] 
In this case the matrix $b(\Gamma)/r(\Gamma)$ acts on the right on row vectors.
\end{definition}

The $p$-adic condition that we shall need is:\\

{\it (O.C.) The connections $\nabla$ and $\check{\nabla}$ are overconvergent.}\\

\begin{theorem}[Baldassarri-Chiarellotto]\label{Thm-CT}
Let the pairs $(\E,\nabla)$ and $(\E^\dagger,\nabla^\dagger)$ be defined
as in Section \ref{Sec-DefEetc}. Assume that
conditions (Rat.) and (O.C.) are met.  Then the natural
morphism $(\E,\nabla) \rightarrow (\E^\dagger,\nabla^\dagger)$ induces
an isomorphism $E^{1,n}_{2,dR} \cong E^{1,n}_{2,rig}$.
\end{theorem}

\begin{proof}
This is an application of \cite[Corollary 2.6]{BC}.
\end{proof}

We shall give an effective/algorithm proof of Theorem \ref{Thm-CT} under
the stronger assumption (Prep.Rat.). More precisely, in Theorem \ref{Thm-Growth}
we give effective
bounds on the $p$-adic growth of forms during the reduction algorithm in
the proof of Theorem \ref{Thm-FGEdR} under assumptions
(Prep.Rat.) and (O.C.). It is easy to deduce surjectivity
of the morphism $E^{1,n}_{2,dR} \rightarrow E^{1,n}_{2,rig}$
from Theorem \ref{Thm-Growth}; injectivity may also be easily derived
using the technique in the proof of \cite[Theorem 2]{AA}. We omit
the details of the proof of the isomorphism $E^{1,n}_{2,dR} \cong E^{1,n}_{2,rig}$ from
Theorem \ref{Thm-Growth} since these are not useful to us.

We introduce notation needed for the statement of Theorem \ref{Thm-Growth}:
For each $k$ with $1 \leq k \leq m$ denote by $e_k$ the element of
$K^m$ with $1$ in positiion $k$ and $0$ elsewhere.
For $\ell \geq \rho$, $0 \leq j < d$ and $1 \leq k \leq m$ define
\[ u^{(r,\ell,j,k)}  = e_k \frac{\Gamma^j}{r^{\ell+1}} \otimes d\Gamma \in \E \otimes d\Gamma.\]
Apply the algorithm in the first stage of the proof of Theorem
\ref{Thm-FGEdR} to compute 
\[ v^{(r,\ell,j,k)} = \sum_{i = \rho}^{\ell} \frac{V^{(r,\ell,j,k)}_i}{r^i},\, V^{(r,\ell,j,k)}_i 
\in K[\Gamma]^m, \,\deg(V^{(r,\ell,j,k)}_i) < d\]
such that $u^{(r,\ell,j,k)} - \nabla(v^{(r,\ell,j,k)}) =: w^{(r,\ell,j,k)} 
\in \E \otimes d\Gamma_\rho$. Similarly, for any $\ell \geq \rho$ we may apply
the algorithm in the second stage of the proof of Theorem \ref{Thm-FGEdR} to write
\[  u^{(\infty,\ell,k)}  - \nabla(v^{(\infty,\ell,k)}) = w^{(\infty,\ell,k)} \in \E \otimes d\Gamma_\rho\]
where this time
\[ u^{(\infty,\ell,k)} := e_k \Gamma^{\ell - 1} \otimes d\Gamma,\mbox{ and }
v^{(\infty,\ell,k)} = \sum_{i = \rho}^\ell V^{(\infty,\ell,k)}_i \Gamma^i \mbox{ for some } V^{(\infty,\ell,k)}_i \in K^m.\]
We have the following effective bounds on the growth of forms during reduction c.f.
\cite[Lemma 2]{KK1}.\footnote{For notational convenience, in the statement of the theorem
and also inequalities (\ref{Eqn-Target}), (\ref{Eqn-Local}), and (\ref{Eqn-DetBound}) the expression
$\log_p(\ell - \rho)$ occurs, when the argument could be zero; similarly $\log_p(j)$ where
$j = 0$ occurs in inequality (\ref{Eqn-ExpInt}). In these cases ``$\log_p(0)$'' should be understood
to be zero.}

\begin{theorem}\label{Thm-Growth}
Assume that conditions (Prep.Rat.) and (O.C.) hold.
Then for $w \in \{w^{(r,\ell,j,k)},w^{(\infty,\ell,k)}\}$ we have
\[
\ord_p(w) \geq - \left( \alpha \log_p(\ell - \rho) + \beta \right)
\]
for some effective constants $\alpha, \beta \in \Q$ which depend only upon the connection $\nabla$ and
the basis $\B$, i.e., are independent of  the starting form $u  \in \{u^{(r,\ell,j,k)},u^{(\infty,\ell,k)}\}$.
\end{theorem}

We shall make the constants $\alpha,\beta$ completely explicit
in Note \ref{Note-Constants}.
We note that Theorem \ref{Thm-Growth} also holds with the same constants if one applies the variant algorithm for reducing
pole orders at infinity given at the end of the proof of Theorem \ref{Thm-FGEdR}.
Since the forms $u^{(r,\ell,j,k)}$ and $u^{(\infty,\ell,k)}$ span $\E \otimes d\Gamma$ as a $K$-vector
space, the above theorem allows one to deduce bounds on the growth of arbitrary
forms during the reduction algorithm.

The proof of Theorem \ref{Thm-Growth} will be reduced by a localisation argument
to that of giving effective bounds on the $p$-adic convergence of the uniform part of the
local solution matrix to a differential system at a regular singular point. Such bounds are
provided in Lemma \ref{Lem-DGS}, whose proof in turn relies on a deep theorem of
Christol-Dwork-Gerotto-Sullivan \cite[Chap V]{DGS}, and an elementary result of Clark (Lemma \ref{Lem-Naive}). 

\begin{proof}
Since $b(\Gamma) \in M_m (\Oh_K[\Gamma])$ and $\ord_p(r(\Gamma)) = 0$,
from the equation ``$w = u - \nabla(v)$'' we see that it suffices to prove that for
$v \in \{v^{(r,\ell,j,k)},v^{(\infty,\ell,k)}\}$ we have
\begin{equation}\label{Eqn-Target}
\ord_p(v) \geq - \left( \alpha \log_p(\ell - \rho) + \beta \right)
\end{equation}
for some effective constants $\alpha, \beta \in \Q$ which depend only upon the connection $\nabla$ and
the basis $\B$.

We divide the proof of (\ref{Eqn-Target}) into three steps:
\begin{itemize}
\item{{\bf Step 1}: We reduce proving bound (\ref{Eqn-Target}) to proving
the local bounds (\ref{Eqn-Local}). Here we need that $r$ is
squarefree modulo $p$.}
\item{{\bf Step 2}: We reduce proving each local bound (\ref{Eqn-Local}) to proving
a different local bound (\ref{Eqn-UniformBound}). Here we need assumption (Prep.Rat.).}
\item{{\bf Step 3}: Bound (\ref{Eqn-UniformBound}) is deduced from an effective version of a theorem
of Christol. This step uses assumptions (Prep.Rat.) and (O.C.).}
\end{itemize}

Recall that $\rho$ is defined to be the smallest
integer larger than all integers in the exponent set
$E(\nabla,\B)$. We note that the argument we give works for {\it any} ``$\rho$'' larger than every integer
in the exponent set $E(\nabla,\B)$.

{\bf Step 1}: First let us consider the case that $u := u^{(r,\ell,j,k)}$ for some $\ell \geq \rho$,
$0 \leq j < d$ and $1 \leq k \leq m$. 
Let us simplify the notation above by removing the exponent
``$(r,\ell,j,k)$'' where it occurs, i.e., $v := v^{(r,\ell,j,k)}$, $w := w^{(r,\ell,j,k)}$ and $V_i := V_i^{(r,\ell,j,k)}$ etc. 
Let $\gamma \in R$ be a root of $r(\Gamma) = 0$. Let $t_\gamma = \Gamma - \gamma$
and expand $v$ locally as
\begin{equation}\label{Eqn-Localv}
 v(t_\gamma) = v_{\gamma,\ell} t_\gamma^{-\ell} + \dots + v_{\gamma,\rho} t_\gamma^{-\rho} + \dots,\, 
v_{\gamma,i} \in K(\gamma)^m.
\end{equation}
We show now that (\ref{Eqn-Target}) holds for $v = v^{(r,\ell,j,k)}$ provided that 
\begin{equation}\label{Eqn-Local}
 \ord_p(v_{\gamma,i}) \geq -(\alpha \log_p(\ell - \rho) + \beta) \mbox{ for all $\gamma \in R$ and
$\rho \leq i \leq \ell$.}
\end{equation}
Assume (\ref{Eqn-Local}) holds. We claim that $\ord_p(V_i) \geq -(\alpha \log_p(\ell - \rho) 
+ \beta)$ for $\rho \leq i \leq \ell$, from which (\ref{Eqn-Target}) follows immediately. This claim can be proved by descending induction on
$i$ in this range. Less formally, observe that $v_{\gamma,\ell}$ is just
$V_\ell(\gamma) r^\prime(\gamma)^{-1}$.
Since the roots of $r$ are distinct modulo $p$ we have $\ord_p(V_\ell(\gamma)) = 
\ord_p(v_{\gamma,\ell})$. Since $\deg(V_\ell) < d = |R|$, from (\ref{Eqn-Local}) we deduce the claimed
bound on $\ord_p(V_\ell(\Gamma))$. Now subtract $V_\ell/r^\ell$ from both sides
of (\ref{Eqn-Localv}) and repeat the argument for $i = \ell  - 1$, and so on.

Similarly, assuming (\ref{Eqn-Local}) holds for the coefficients in the local expansion of $v:= v^{(\infty,\ell,k)}$
at infinity, we easily deduce that (\ref{Eqn-Target}) holds for $v = v^{(\infty,\ell,k)}$. 

It remains to establish the local bound (\ref{Eqn-Local}). (We omit the remainder
of the proof for $v = v^{(\infty,\ell,k)}$ since it is exactly the same.)

{\bf Step 2}: Fix $\gamma \in R$ and simplify notation as in Step 1.  
Define $G(t) \in M_m(\Oh_{K(\gamma)})[[t]]$ so that $-t^{-1}G(t)$ is the 
expansion of $b(\Gamma)/r(\Gamma)$  w.r.t the local parameter $t := t_\gamma = \Gamma - \gamma$.
Define $H := G(0)$, the negative of the residue matrix $b(\gamma)/r^\prime(\gamma)$.
Let the local solution matrix $Y(t) t^H$ to the differential system
$t \frac{d}{dt} - G(t) = 0$ be defined as in \cite[Chap. III Prop. 8.5]{DGS}.
(Note that we have chosen our signs to be consistent with \cite{DGS}.)
The existence of such a solution matrix requires the assumption
(Prep.Rat.). The {\rm uniform
part} $Y(t)$ lies in $M_m(K(\gamma))[[t]]$ with $Y(0) = I_m$, and the element $t^H$, which
is constructed on \cite[Page 103]{DGS}, satisfies the equation
$\frac{d}{dt}(t^H) = t^{-1}t^H$, and $(t^H)^{-1} = t^{-H}$.

Write $Y = \sum_{i = 0}^\infty Y_i t^i$ and
$Y(t)^{-1} = \sum_{i = 0}^\infty Z_i t^i$, where $Y_i,Z_i \in M_m(K(\gamma))$. We shall
show now that it is enough to prove that
\begin{equation}\label{Eqn-UniformBound}
\ord_p(Y_i),\,\ord_p(Z_i) \geq -(\alpha_1 \log_p(i) + \beta_1),\,i \geq 1
\end{equation}
for some explicit $\alpha_1,\beta_1 \in \Q$ which
depend only on $\nabla$ and $\B$.

Let us assume (\ref{Eqn-UniformBound}) holds.  
Observe that we have the local factorisation
\begin{equation}\label{Eqn-FactorisationSP}
 \nabla = Y(t) t^H \circ \frac{d}{dt} \circ t^{-H} Y(t)^{-1}.
 \end{equation}
Premultiplying the localised equation $\nabla(v) = u-w$ by $(Y(t) t^H)^{-1}$, using 
(\ref{Eqn-FactorisationSP}), and then integrating,
we find that
\begin{equation}\label{Eqn-Integral}
 v(t) = Y(t) t^{H} \left\{ \int t^{-H}  Y(t)^{-1} (u - w) dt + c \right\},
 \end{equation}
 for some constant $c \in K(\gamma)^m$.
Bound (\ref{Eqn-Local}) can now be deduce by explicitly integrating the
righthand side of (\ref{Eqn-Integral})
and comparing coefficients of $t^{-i}$ for $\rho \leq i \leq \ell$.

Specifically,  the integrand on the righthand side of (\ref{Eqn-Integral}) can be written 
\begin{equation}\label{Eqn-ExpInt}
 \sum_{j \geq  0} t^{-H} a_j t^{-(\ell+1) + j},\, \ord_p(a_j) \geq -(\alpha_1 \log_p(j) + \beta_1) \mbox{ for }
0 \leq j < \ell+1 - \rho,
\end{equation}
for some $a_j \in K(\gamma)^m$. The lower bound on $\ord_p(a_j)$ comes
from (\ref{Eqn-UniformBound}) and the integrality of $u(t)$.
Note that we do not have any bounds on $\ord_p(a_j)$ for $j \geq \ell +1 - \rho$ since
these terms are affected by the unknown element $w(t)$ and unknown constant $c$.  
Element (\ref{Eqn-ExpInt}) may be
explicitly integrated ``term--by-term''. Precisely, from the defining property of the element
$t^H$, ones sees that $\int t^{-H} a_j t^{-(\ell + 1) + j} dt = (-H - (\ell + 1) + j + 1)^{-1}
t^{-H} a_j t^{-(\ell + 1) + j + 1}$, plus an unknown constant of integration. Recall
that $-H$ is the residue matrix. Now
for $0 \leq j < \ell + 1 - \rho$ we have that $-\ell \leq -(\ell + 1) + j + 1 \leq -\rho$, and since
$\rho$ is larger than any eigenvalue of $-H$ the inverse matrix immediately
above exists.
 (For $j \geq \ell + 1 - \rho$, when the inverse does not exist the coefficient
$a_j$ must be zero.) Next, note that $(-H+i)^{-1}$ for $i \in \Z$ (when it exists)
commutes with $t^H$; this follows from the fact that $H$ commutes
with $t^H$, see \cite[Page 103]{DGS}. Thus each term on the righthand side of
(\ref{Eqn-Integral}) which does not involve the constant $c$ has the form
\begin{equation}\label{Eqn-Mess}
Y_i t^i (t^H) (-H -(\ell + 1)  + j + 1)^{-1} t^{-H} a_j t^{-(\ell + 1) + j + 1}
 = Y_i (-H - \ell + j)^{-1} a_j t^{-\ell + i  +  j}
 \end{equation}
for some $i,j \geq 0$. Terms on the righthand side of (\ref{Eqn-Integral})
which do involve $c$ have the form $Y_i t^i t^{H} c$ for some $i \geq 0$.
Since the lefthand side $v(t)$ is a Laurent series, it follows from
\cite[Chap V Lemma 2.3]{DGS} that either $c = 0$, or $c \ne 0$ with $H$ a diagonal
matrix with integer eigenvalues. Since all eigenvalues of $-H$ are less than $\rho$, 
in either case any term on the righthand side of (\ref{Eqn-Integral})
involving $c$ cannot effect the coefficient
of $t^{-i}$ for $\rho \leq i \leq \ell$.

From (\ref{Eqn-Mess}), a lower bound on the coefficient of $t^{-\ell  + s}$  for
$0 \leq s  \leq \ell  - \rho$  on the righthand side of (\ref{Eqn-Integral}) is
\begin{equation}\label{Eqn-Min}
 \min_{i + j = s } \ord_p \left( Y_i (-H - \ell + j)^{-1} a_j \right).
 \end{equation}
We have bounds on $\ord_p(Y_i)$ and $\ord_p(a_j)$ for $i \geq 0$ and
$0 \leq j \leq \ell - \rho$, viz. (\ref{Eqn-UniformBound}) and (\ref{Eqn-ExpInt}).
It remains to bound $\ord_p((-H -\ell + j)^{-1})$ for $0 \leq j \leq \ell - \rho$.
Now $\ord_p((-H - \ell + j)^{-1}) \geq - \ord_p(\det(-H - \ell + j))$ so we must
find an upper bound for the valuation of the determinant. 
Denote by $\lambda_1,\dots,\lambda_m \in \Q$ the eigenvalues of
$-H$ (the residue matrix). Then $\det(-H - \ell + j) = \prod_{i = 1}^m (\lambda_i - \ell + j)$. 
Take the positive integer $N$ to be a lowest common denominator for the $\lambda_i$ and define
$\mu_i = N\lambda_i \in \Z$; note
that $\gcd(p,N) = 1$ since the eigenvalues are $p$-adic integers. Take
the positive integer $\Delta$ to be minimal so that $|\lambda_i| \leq \Delta$
for all $i$. Then for $0 \leq j \leq \ell - \rho$ we have
\[ \ord_p(\lambda_i - \ell + j) = \ord_p(\mu_i - N(\ell - j)) \leq
\log_p( |\mu_i| + N\ell)  \leq \log_p(N) + \log_p(\Delta + \ell).\]
So certainly for $0 \leq j \leq \ell - \rho$ we have
\begin{equation}\label{Eqn-DetBound}
\ord_p(\det(-H - \ell + j)) \leq m \left( \log_p(\ell - \rho) + \log_p(N) + \log_p(2\Delta+2) \right),
\end{equation}
since $\rho \leq \Delta + 1$.
From (\ref{Eqn-Min}) and (\ref{Eqn-DetBound}) we conclude that
\begin{equation}\label{Eqn-abab1}
 \alpha := \lceil 2\alpha_1 + m \rceil,\, \beta := \lceil 2\beta_1 + m(\log_p(N) + \log_p(2\Delta +2))
\rceil
\end{equation}
will certainly suffice.

{\bf Step 3}: We now establish bound (\ref{Eqn-UniformBound}). 
First consider $Y(t)$. By assumptions (Prep.Rat.) and
(O.C.), we see that the hypothesis of Lemma \ref{Lem-DGS}
are met. Hence we may apply the bound in Lemma \ref{Lem-DGS} to our
differential system $t \frac{d}{dt} - G(t) = 0$.
We next note that $Y(t)^{-1}$ is the uniform part of
the local solution matrix of the ``dual'' differential system
$t\frac{d}{dt} + G(t) = 0$ where $G(t)$ acts on the right, or equivalently, the
transpose of the uniform part of the local solution matrix of $t\frac{d}{dt} + G(t)^{tr} = 0$ with
$G(t)^{tr}$ acting on the left, c.f. \cite[Page 193]{DGS}. 
Again by assumptions (Prep.Rat.) and
(O.C.) the hypothesis of Lemma \ref{Lem-DGS} are met,
and we may use that bound.
\end{proof}

The next lemma is a modest generalisation of the main theorem in
\cite[Chap. V]{DGS}. 

\begin{lemma}\label{Lem-DGS}
Let $G(t) \in M_m(\Oh_{K(\gamma)})[[t]]$ be the local expansion of a rational function
$G(\Gamma) \in M_m(K(\Gamma))$ around some point $t = \Gamma - \gamma$.
Let $\delta := t \frac{d}{dt}$ and consider the differential system
$\delta - G = 0$. Assume
\begin{enumerate}
\item{The eigenvalues of $G(0)$ are rational numbers, and no two differ
by a positive integer,}
\item{The solution matrix to the differential system around the generic
point $t$ with $|t|_p = 1$ converges $p$-adically on the open
unit disk.}
\end{enumerate}
 Let $Y(t) = 
\sum_{i = 0}^\infty Y^i t^i \in M_m(K(\gamma))[[t]]$ be the uniform
part of the local solution matrix $Y(t) t^{G(0)}$ of the differential system
$\delta - G = 0$.
Then there exist $\alpha_1,\beta_1 \in \R$ such that
\[ \ord_p(Y_i) \geq -(\alpha_1 \log_p(i) + \beta_1),\mbox{ for all }\,i \geq 1. \]
\end{lemma}

\begin{proof}
First, change basis by a matrix
$\h \in GL_m(K(\gamma)[t,t^{-1}])$ so that 
\[ G_{[\h]} := \h^{-1} t\frac{d\h}{dt} + \h^{-1} G\h\]
is such that $G_{[\h]}(0)$ has eigenvalues
in the interval $[0,1)$. By \cite[Chap. V, Prop. 4.1]{DGS} the matrix
$\h$ may be taken to be unimodular, i.e., $\ord_p(\h),\ord_p(\h^{-1}) \geq 0$.
Moreover, the degree in $t$ of $\h$ (degree in $t^{-1}$ of $\h^{-1}$) is the
absolute value of the floor of the most negative eigenvalue of $G(0)$, and the
degree in $t^{-1}$ of $\h$ (degree in $t$ of $\h^{-1}$) is the floor of the most
positive eigenvalue of $G(0)$; this is easily seen by viewing $\h$ and $\h^{-1}$
as a product of shearing transformations.

Let $\tilde{Y}(t) = \sum_{i = 0}^\infty \tilde{Y}_i t^i$ be the uniform
part of the local solution matrix to the differential
system $\delta - G_{[\h]}$. Note that $G_{[\h]} \in M_m(\Oh_{K(\gamma)})[[t]]$ by
the unimodularity of $\h$.

The main theorem in
\cite[Chap. V Section 9]{DGS} assures us that there
exists $\alpha_2, \beta_2 \in \R$ such that
\[ \ord_p(\tilde{Y}_i) \geq -(\alpha_2 \log_p(i) + \beta_2),\, i \geq 1.\]
We comment briefly on why the main theorem is applicable:
In the notation of \cite[Chap V]{DGS}, we must check conditions $\cR 1$,
$\cR 2$, $\cR 3^\prime$ and $\cR 4$. Now $\cR 1$ is true since the matrix
$G_{[\h]}$ contains functions which are localisations of rational functions; $\cR 2$ 
(overconvergence) follows from assumption 2. and the unimodularity of $\h$;
$\cR 2$ (eigenvalues in $\Z_p \cap [0,1)$) is true by assumption 1.; $\cR 4$ (integrality)
follows since $G_{[\h]} \in M_m(\Oh_{K(\gamma)})$.

There exists
$\bar{\h} \in GL_m(K[t,t^{-1}])$ such that $Y(t)
= \h^{-1} \tilde{Y} \bar{\h}$, c.f. \cite[Page 163 Lines 1-6]{DGS}.
Moreover, the degree in $t$ of $\bar{\h}$ is the floor
of the most positive eigenvalue of $G(0)$, and the degree in $t^{-1}$
of $\bar{\h}$ is the absolute value of the floor of the most negative
eigenvalue of $G(0)$; this follows from the argument on
\cite[Page 163 Lines 11-21]{DGS}.

Since
$\ord_p(\h^{-1}) \geq 0$, to prove the lemma we need only calculate a lower bound
on $\ord_p(\bar{\h})$. One computes this from the equation
$\bar{\h} = \tilde{Y}^{-1} \h Y$ using bounds on the degree
in $t^{-1}$ of $\h$ and the naive upper bound on the growth of the coefficients of
$\tilde{Y}(t)^{-1}$ and $Y(t)$, c.f. \cite[Page 191-193]{DGS}.  The naive upper
bounds are given in Lemma \ref{Lem-Naive}. Specifically, one finds certainly
\begin{equation}\label{Eqn-b3}
 \ord_p(\bar{\h}) \geq -\beta_3,\, \beta_3 := m^2 \left(\frac{\Delta}{p-1} + 4 \log_p(\Delta + 1) + 
2\log_p(2N)
\right)
 \end{equation}
where $\Delta$ is such that all eigenvalues $\lambda$ of $G(0)$ have $|\lambda| \leq \Delta$ and
$N$ the lowest common denominator for the eigenvalues. Note that we have already
observed $\deg_{t^{-1}}(\h^{-1}),\deg_{t^{-1}} (\bar{\h}) \leq \Delta$.
Comparing coefficients in
$Y = \h^{-1} \tilde{Y} \bar{\h}$ we see
\[ \ord_p(Y_i) \geq \ord_p(\tilde{Y}_{i+2\Delta}) - \beta_3  \geq 
- (\alpha_2 \log_p(i) + \alpha_2 \log_p(2\Delta + 1) + \beta_2 + \beta_3) \]
which gives a bound of the required form. Precisely, take
\begin{equation}\label{Eqn-ab1ab2}
\alpha_1 := \alpha_2,\,\beta_1 :=  \alpha_2 \log_p(2\Delta + 1) + \beta_2 + \beta_3.
\end{equation}
\end{proof}

The following lemma is an effective version of the general bound of
Clark \cite{DNC,MS}.

\begin{lemma}
\label{Lem-Naive}
Situation as in the statement of Lemma \ref{Lem-DGS}, only without assumption 2. Denote
by $\lambda_1,\dots,\lambda_m$ the eigenvalues of $G(0)$ and assume that 
$|\lambda_i| \leq \Delta$ and $N \lambda_i \in \Z$ 
for all $i$ and some minimal integers $\Delta \geq 0$, $N \geq 1$ with
$\gcd(p,N) = 1$. 
Then for all $s \geq 1$ we have
\[ \ord_p(Y_s) \geq -m^2 \left( \frac{s}{p-1} + \log_p(1+s) + \log_p(2\Delta + 1) + \log_p(N)\right).\]
\end{lemma}

\begin{proof}
The power series $Y(t)$ may be computed using the classical method
in \cite[Chap V, Remark 2.2]{DGS}. It follows from this recursive method that
\[ \ord_p(Y_s) \geq - \sum_{k = 1}^s \sum_{i,j = 1}^m \ord_p(k - \delta_{ij})\]
where $\delta_{ij} := \lambda_i - \lambda_j$. We shall estimate the righthand side using
the method of Clark. Specifically, for $\alpha \in \Zp$ and $s$ a positive integer,
define $\theta(\alpha,s):=\prod_{k = 1}^s (\alpha + k)$. Then $\ord_p(Y_s) \geq - \sum_{i,j}
\ord_p\theta(-\delta_{ij},s)$. We compute an upper bound for the $\theta(-\delta_{ij},s)$ using
the argument in \cite[Page 265, Case 3]{DNC}. First, write $-\delta_{ij} = \nu_{ij}/N$, where
$|\nu_{ij}| \leq 2N \Delta$. Then for any positive integer $s$ we have 
\[ \ord_p(-\delta_{ij} + s) = \ord_p(\nu_{ij} + Ns) \leq \log_p(|\nu_{ij}| + Ns) \leq
\log_p(N) + \log_p(2\Delta+1)  + \log_p(s).\]
This shows that we may take the expression ``$v(x) = -k\log_p(1 + x) - k^\prime$'' 
immediately preceding \cite[Eqn (13)]{DNC}, to have coefficients ``$k = 1$'' and
``$k^{\prime} = \log_p(2 \Delta + 1) + \log_p(N)$'' (we have changed Clark's
``$\log$'' to ``$\log_p$''). From \cite[Eqn (14)]{DNC} we deduce that
\[ \ord_p \theta(-\delta_{ij},s) \leq \frac{s}{p-1} - v(s)  = \frac{s}{p-1} +\log_p(1 + s) +
\log_p(2\Delta + 1) + \log_p(N)\]
as required.
\end{proof}

\begin{note}\label{Note-Constants}
The constants $\alpha,\beta \in \Q$ in Theorem \ref{Thm-Growth} can be made completely
explicit. Precisely, by equations (\ref{Eqn-abab1}), (\ref{Eqn-b3}) and (\ref{Eqn-ab1ab2}) one sees that it suffices to make the constants $\alpha_2$ and $\beta_2$ from the proof of Lemma
\ref{Lem-DGS} explicit. These constants are those which occur in the theorem of
Christol-Dwork-Gerotto-Sullivan.

The theorem of CDGS states that
$\ord_p(Y_i) \geq -(\alpha_2 \lfloor \log_p(i) \rfloor + \beta_2)$ with $\alpha_2,\beta_2
\in \R$ as follows. Define
\begin{equation}\label{Eqn-Bmp}
 B_{m,p}:= m-1+\ord_p((m-1)!) + \min\left\{m-1,\ord_p \prod_{j = 1}^m {m \choose j}\right\}.
 \end{equation}
When all the eigenvalues are zero (nilpotent case) we can directly apply the
Christol-Dwork theorem \cite[Chap. V, Thm 2.1]{DGS}. This gives
$\alpha_2 = B_{m,p}$ and $\beta_2 =0$.
In particular, for $p \geq m$ we have $\alpha_2 = m-1$, and in general $\alpha_2 \leq 3m-3$.
For eigenvalues in the interval $[0,1)$, one applies the generalisation
in \cite[Sec. 9, Chap. V]{DGS}. Define $\ell := \lfloor \log_p(i) \rfloor + 1$. From the equation on middle
of page 198 in \cite{DGS} one deduces that
\[ \ord_p(Y_i) \geq -\left\{(\ell + 1)(m(m-1)) + \ell B_{m,p} + B\right\} \]
where the number $B$ is defined on  \cite[Pages 197-198]{DGS}.
(Note that we have changed from multiplicative to additive notation, so
our ``$B$'' corresponds to ``$\log_p(B)$''.) Let $N$ be the lowest  common denominator
for the eigenvalues of the residue matrix. Then for $p > \{m,2N\}$ from \cite[Remark 2.2]{DGS} one sees that $B_{m,p} = m-1$ and $B = 0$, so 
$\alpha_2 = m^2 - 1$ and $\beta_2 = 2m^2 - m - 1$.
For general $p$, one computes from \cite[Page 198, Line 7]{DGS} and
\cite[Page 197, Line 12]{DGS}
 that $B \leq (\ell+1)m(m-1)\log_p(2N)$, and as observed before
$B_{m,p} \leq 3m-3$. Thus we may take
\[ \alpha_2 = m(m-1)(1 + \log_p(2N)) + 3m-3,\, \beta_2 = 2m(m-1) ( 1+ \log_p(2N)) + 3m-3.\]
We conclude the following: if $\Delta \geq 0$ is a bound on the absolute value of the local monodromy eigenvalues,
$N \geq 1$ a lowest common denominator, and $m$ the dimension, then one has
\begin{equation}\label{Eqn-AbsEst}
 \alpha = C_1m^2 (1 + \log_p(N)), \, \beta = C_2 m^2 \left(1 + \frac{\Delta}{p} + 
 \log_p (\Delta+1)  + \log_p(N)\right)
 \end{equation}
for some {\it absolute} constants $C_1,C_2 \in \Q$. This will be useful in
our complexity estimates.
\end{note}

In Section \ref{Sec-CTholds} we shall see that conditions (Rat.) and (O.C.) are met
in the situations (which arise ``from geometry'') which we shall encounter. The stronger
condition (Prep.Rat.) will be met in the examples we compute in Section \ref{Sec-CompData}.

Theorem \ref{Thm-Growth} is essential in the complexity analysis and practical application of the algorithm in
the proof of Theorem \ref{Thm-FGEdR} for the following reason: When one calculates the reduction
of differential forms using this algorithm it is impractical to store the coefficients ``exactly''. At each
step of the reduction, one ``approximates'' the coefficients modulo some fixed power of the
characteristic. Making this approximation amounts to adding an ``error form'' to the form being reduced.
Theorem \ref{Thm-Growth} shows that the error introduced propagates in a ``logarithmic'' manner
during the remainder of the computation. Furthermore, Theorem \ref{Thm-Growth} applied
with a general value ``$\rho$'',  at least as big as $\rho$ itself, allows one to bound the intermediate 
coefficient size during the reduction computation. Note that a naive inspection of the reduction
formulae in the proof of Theorem \ref{Thm-FGEdR} suggests that terms grow and errors
propagate in a ``linear'' manner; that this is {\it not} the case for calculations in rigid cohomology was an important insight of Kedlaya c.f. \cite[Lemma 2]{KK1}.

A similar ``logarithmic error propogation'' phenomenon arises in the numerical
solution of differential systems, as we shall see in Theorem \ref{Thm-DiffSysLoss}
of the next section.

\section{Deformation of Frobenius}\label{Sec-DefFrob}

In this section we retain the notation and definitions in the first paragraph of Section \ref{Sec-DefEetc},
but slightly alter some of our assumptions.
Specifically, $\E$ is a free $A$-module of rank $m$, where $A = K[\Gamma,1/r(\Gamma)]$ and
$K$ is the unramified extension of $\Qp$ of degree $[\fq:\fp]$. The map $\nabla:
\E \rightarrow \E \otimes \Omega^1_A$ is a connection which with respect to a fixed
basis of $\E$ (and ``natural'' corresponding basis of $\E \otimes \Omega^1_A)$ acts as:
\[ \nabla = \frac{d}{d\Gamma} + B(\Gamma),\, B(\Gamma) = \frac{b(\Gamma)}{r(\Gamma)}.\]
Here $b(\Gamma) \in M_m(\Oh_K[\Gamma])$. In this section we {\it do not} assume that the
connection has only simple poles, i.e., we do not need the assumption that $r(\Gamma) \in
\Oh_K[\Gamma]$ is squarefree, nor do we need the degree restriction on the matrix $b(\Gamma)$.
However, we shall add the new assumption that $r(0) \not \equiv 0 \bmod{p}$.

Let $\nabla^\dagger: \E^\dagger 
\rightarrow \E^\dagger \otimes \Omega^1_{A^\dagger}$ be obtained from $\nabla$ by
extension of scalars; the ring $A^\dagger = K[\Gamma,1/r(\Gamma)]^\dagger$ is described explicitly
in Section \ref{Sec-Pencils}. Let $\sigma: A^\dagger \rightarrow A^\dagger$ be the lifting of
the $p$th power Frobenius which maps $\Gamma \mapsto \Gamma^p$.
Let $F: \E^\dagger \rightarrow \E^\dagger$ be an injective $\sigma$-linear map such that the triple
$(\E^\dagger,\nabla^\dagger,F)$ defines a $(\sigma,\nabla)$-module over $A^\dagger$ (see Section 
\ref{Sec-CMD}), i.e., the diagram (\ref{Diagram-F}) commutes.

\subsection{Local deformation}\label{Sec-LocalDef}

If we assume that the connection matrix $B(\Gamma)$ is known and that the specialisation
$F(0)$ is also known, the commutative diagram (\ref{Diagram-F}) allows
the computation of a {\it local} expansion of the Frobenius matrix $F(\Gamma)$ around the
origin to any required precision. We describe two different approaches.

\subsubsection{Method 1}\label{Sec-Method1}

Let $C(\Gamma)$ be a basis of local solutions to the differential system $\nabla = 0$ with
initial condition $C(0) = I_m$. So
\begin{equation}\label{Eqn-C}
\frac{dC}{d\Gamma} + B(\Gamma)C(\Gamma) = 0.
\end{equation}
Commutativity of (\ref{Diagram-F}) implies that the Frobenius map $F$ is stable on the kernel of the connection. Recalling that the map $F$ is $\sigma$-linear, we deduce the matrix equation 
\[ F(\Gamma)C^\sigma(\Gamma^p) = C(\Gamma)D \]
where $D$ is some constant matrix. Evaluating both sides at $\Gamma = 0$ shows $D = F(0)$. So we have the local factorisation 
\begin{equation}\label{Eqn-F}
 F(\Gamma) = C(\Gamma) F(0) (C^\sigma(\Gamma^p))^{-1}.
 \end{equation}
Thus $F(\Gamma)$ can be computed modulo $\Gamma^{N_\Gamma}$ for any $N_\Gamma \geq 1$ 
provided we can compute
$C(\Gamma)$ modulo $\Gamma^{N_\Gamma}$. A simple recursion formula for computing the matrix coefficients in the 
local expansion of $C(\Gamma) = \sum_{\ell = 0}^\infty C_\ell \Gamma^\ell$ can be derived from the equation
\[ r(\Gamma) \frac{dC}{d\Gamma} + b(\Gamma) C(\Gamma) = 0.\]
Specifically, write $r(\Gamma) = \sum_{i = 0}^{\deg(r)} r_i \Gamma^i$ and
$b(\Gamma) = \sum_{i = 0}^{\deg(b)} b_i \Gamma^i$. Then for  $\ell \geq 1$ we have 
\begin{equation}\label{Eqn-Recursion}
C_\ell = -\frac{1}{r_0 \ell} \left(\sum_{i= 0}^{\deg(b)}  b_i C_{(\ell - 1)-i} + \sum_{i = 1}^{\deg(r)} r_i
(\ell - i) C_{\ell - i}\right).
\end{equation}
One can, of course, compute the series $C^\sigma(\Gamma^p)^{-1}$ by power series inversion. However,
it is better to observe that the matrix $C(\Gamma)^{-1}$ is the solution of the
``dual equation"
\begin{equation}\label{Eqn-DualC}
 \frac{d C^{-1}}{d\Gamma} - C(\Gamma)^{-1}B(\Gamma) = 0,\, C(0)^{-1} = I_m.
 \end{equation}

It is impractical to carry out the above computations using ``exact arithmetic''; one desires to
``truncate'' each coefficient $C_\ell$ ``modulo $p^N$'' for some appropriate $N > 0$ after it has
been computed. It is an essential task to analyse the ``propagation'' of the error this
introduces as one continues the computation. 

Let $E_\ell \in M_m(\Oh_K)$ for $\ell \geq 1$, and $N$ be a non-negative integer. 
Let the sequence $D_\ell \in M_m(K)$ for
$\ell \geq 0$ be computed in the following manner. Define $D_0 := I_m$ and
for $\ell \geq 1$,
\begin{equation}\label{Eqn-ApproxSol}
 D_\ell := -\frac{1}{r_0 \ell} \left(\sum_{i= 0}^{\deg(b)}  b_i D_{(\ell - 1)-i} + \sum_{i = 1}^{\deg(r)} r_i
(\ell - i) D_{\ell - i}\right) + p^N E_\ell. 
\end{equation}
The series $D(\Gamma) :=\sum_{\ell = 0}^\infty D_\ell \Gamma^\ell$ is an ``approximate
solution'' to the differential system (\ref{Eqn-C}) computed ``modulo $p^N$''. 
In practice, the ``error sequence''
$E_\ell$ is chosen to ensure the $p$-adic expansions of the entries of $D_\ell$ are
``truncated modulo $p^N$''.

\begin{theorem}\label{Thm-DiffSysLoss}
Let $C(\Gamma)$ be the solution to (\ref{Eqn-C}) with $C(0) = I_m$, and $D(\Gamma)$ the
``approximate solution'' defined via equation (\ref{Eqn-ApproxSol}). Then for $\ell \geq 0$ we have
\[ \ord_p(D_\ell - C_\ell) \geq N - \alpha^\prime \lfloor \log_p(\ell+1) \rfloor\]
for some explicitly computable constant $\alpha^\prime \geq 0$. Furthermore, one can take
$\alpha^\prime = 6m-5$ for any prime $p$, and
$\alpha^\prime = 2m - 1$ when $p \geq m$.
\end{theorem}

\begin{proof}
First observe that we have the local factorisation around the origin
\begin{equation}\label{Eqn-LocalFactor}
\nabla = C(\Gamma) \circ \frac{d}{d\Gamma} \circ C(\Gamma)^{-1}.
\end{equation}
Next observe that the series $D(\Gamma)$ is a local solution to the inhomogeneous
differential equation
\[ r(\Gamma) \frac{d D}{d\Gamma} + b(\Gamma)D(\Gamma) = p^N r_0 \sum_{\ell = 1}^\infty \ell E_\ell
 \Gamma^\ell.\]
Thus $\nabla(D) = p^N E(\Gamma) r(\Gamma)^{-1}$ where $E(\Gamma) :=
r_0 \sum_{\ell = 1}^\infty \ell E_\ell \Gamma^\ell$.
Using the local factorisation (\ref{Eqn-LocalFactor})
one deduces that
\[ \frac{d}{d\Gamma} \left( C(\Gamma)^{-1} D(\Gamma)\right) = C(\Gamma)^{-1} p^N E(\Gamma) r(\Gamma)^{-1}.\]
Integrating we find that there exists a constant matrix $c$ such that
\[ D(\Gamma) = C(\Gamma) \left( \int C(\Gamma)^{-1} p^N E(\Gamma) r(\Gamma)^{-1}d\Gamma 
+ c \right)\]
Note that $E(\Gamma) r(\Gamma)^{-1} \in \Gamma \Oh_K[[\Gamma]]$ since
$r(0)$ is a unit.
Since $C(\Gamma) = D(\Gamma) \bmod{\Gamma}$ we deduce that
$c = I_m$. Hence
\[ D(\Gamma) - C(\Gamma) = p^N C(\Gamma) 
\int C(\Gamma)^{-1}  E(\Gamma) r(\Gamma)^{-1} d\Gamma.\]
The connection $\nabla$ and its dual $\check{\nabla}$ come from 
overconvergent $F$-isocrystals, viz, $(\E^\dagger,\nabla^\dagger,F)$ and its ``dual
$(\E^\dagger,\check{\nabla}^\dagger,F^{-1})$''. Hence Dwork's trick of analytic continuation via
Frobenius \cite[Prop. 3.1.2]{NKDwork}
shows that condition (O.C.) is met. Moreover, the connections are regular at
zero, so local exponents are all zero. Thus we
can apply the Christol-Dwork theorem to deduce effective logarithmic bounds on the
growth of coefficients of $C(\Gamma)$ and $C(\Gamma)^{-1}$. Moreover, integration only
has a ``logarithmic'' effect on the growth of coefficients of a power series.  Explicitly, we
can use the constant $B_{m,p}$ in the original Christol-Dwork theorem, see 
(\ref{Eqn-Bmp}) in Note \ref{Note-Constants},
to deduce
$\alpha^\prime = 2B_{m,p} + 1$ and the constant ``$\beta^\prime = 0$''.
This completes the proof.
\end{proof}

Let $N_\Gamma$ and $N$ be  positive integers.
Let $D(\Gamma)$ be an ``approximate solution''  computed ``modulo $p^N$'' to the differential 
system (\ref{Eqn-C}) modulo $\Gamma^{N_{\Gamma}}$. 
Let $\tilde{D}(\Gamma)$ be an ``approximate solution'' computed
``modulo $p^N$'' to the dual system (\ref{Eqn-DualC}) modulo $\Gamma^{\lceil{N_{\Gamma}/p}\rceil}$. Let $G(0) \in M_m(K)$ be such that
$\ord_p(F(0) - G(0)) \geq p^N$. Define $G(\Gamma) := 
D(\Gamma) G(0) \tilde{D}^\sigma(\Gamma^p) \bmod{\Gamma^{N_\Gamma}}$. 
This is our approximation of the local Frobenius
matrix $F(\Gamma)$. We need to bound from below $\ord_p\left((F(\Gamma) \bmod{\Gamma^{N_{\Gamma}}}) - G(\Gamma)\right)$.

From Theorem \ref{Thm-DiffSysLoss}, $D(\Gamma) = C(\Gamma) + p^{N^\prime} e(\Gamma)
\bmod{\Gamma^{N_{\Gamma}}}$ where $\ord_p(e(\Gamma)) \geq 0$ with
$N^\prime := N - (2B_{m,p} + 1) \lfloor \log_p(N_\Gamma) \rfloor$, and
$\tilde{D}^\sigma(\Gamma^p) = (C^\sigma(\Gamma^p))^{-1}  + p^{N^{\prime \prime}} \tilde{e}
(\Gamma) \bmod{\Gamma^{N_{\Gamma}}}$ where $\ord_p(\tilde{e}(\Gamma)) \geq 0$
with $N^{\prime \prime} := N - (2B_{m,p} + 1) \lfloor \log_p(\lceil N_\Gamma/p \rceil) \rfloor$. 
Note that from the Christol-Dwork theorem we have
\[ 
\begin{array}{rcl}
\ord_p\left(C(\Gamma) \bmod{\Gamma^{N_{\Gamma}}}\right) & \geq & -B_{m,p} \lfloor \log_p(N_\Gamma - 1) \rfloor\\
\ord_p\left((C^\sigma(\Gamma^p))^{-1} \bmod{\Gamma^{N_{\Gamma}}}\right) & \geq & -B_{m,p} \lfloor \log_p(\lceil{N_\Gamma/p}\rceil - 1) \rfloor.\\
\end{array}
\]
One now readily calculates a lower bound on
\[ \ord_p\left(\left(C(\Gamma) F(0) (C^\sigma(\Gamma^p))^{-1} \bmod{\Gamma^{N_{\Gamma}}}\right) - 
D(\Gamma) G(0) \tilde{D}^\sigma(\Gamma^p)\right)\]
to be
\[ \min\left\{N^\prime + \ord_p(F(0)) + \ord_p(C^\sigma(\Gamma^p)^{-1} \bmod{\Gamma^{N_\Gamma}}),
\right.\]
\[ N + \ord_p(C(\Gamma) \bmod{\Gamma^{N_\Gamma}}) + 
\ord_p(C^\sigma(\Gamma^p)^{-1} \bmod{\Gamma^{N_\Gamma}}),\]
\[ \left. N^{\prime \prime} + \ord_p(C(\Gamma) \bmod{\Gamma^{N_\Gamma}}) + \ord_p(F(0))\right\}\]
\begin{equation}\label{Eqn-M1loss}
 \geq
N - (3 B_{m,p} + 1) \lfloor \log_p(N_\Gamma) \rfloor + B_{m,p} + \min\left\{\ord_p(F(0)),0\right\}.
\end{equation}
For example, when $F(0)$ has integral entries and $p \geq m$ we have that the
{\it loss} of accuracy when computing $F(\Gamma) \bmod{\Gamma^{N_\Gamma}}$ 
is bounded by $(3m - 2) \lfloor \log_p(N_\Gamma) \rfloor - (m-1)$.

\subsubsection{Method 2}\label{Sec-Method2}

The approach in this section is based upon that taken by Tsuzuki \cite{NT}. We do not give an analysis
of the propagation of errors for this method, although this is an interesting problem.

Commutativity of diagram (\ref{Diagram-F}) implies that 
\[ \frac{dF}{d\Gamma} + B(\Gamma)F(\Gamma) = p\Gamma^{p-1} F(\Gamma)B^\sigma (\Gamma^p).\]
For $\ell \geq 1$, the coefficients $F_\ell$ in the local expansion $F(\Gamma) = \sum_{\ell = 0}^\infty F_\ell \Gamma^\ell$
can be found recursively by rewriting this equation in the form
\[ r(\Gamma)r^\sigma(\Gamma^p) \frac{dF}{d\Gamma} + 
r^\sigma(\Gamma^p)b(\Gamma)F(\Gamma)
 = p \Gamma^{p-1} r(\Gamma) F(\Gamma)b^\sigma (\Gamma^p)\]
and equating the coefficient of $\Gamma^{\ell - 1}$ on both sides. This more direct method
eliminates the multiplication of power series needed to compute the righthand side in (\ref{Eqn-F})
and is also more space efficient.

\subsection{Global deformation: analytic continuation}\label{Sec-GlobalDef}

The entries in the Frobenius matrix $F(\Gamma)$ are $p$-adic holomorphic functions
on the $p$-adic projective line with open unit disks around the poles of $r$ removed, i.e., uniform
limits of rational functions on this closed domain $D_1$, say. Therefore, they can be uniformly approximated on this
domain $D_1$ modulo any power of $p$ by a matrix of rational functions whose denominators
are powers of $r(\Gamma)$. Using the method in Section \ref{Sec-LocalDef}, one can compute
the local expansions of these holomorphic functions to any required $p$-adic and
$\Gamma$-adic accuracy. We now sketch how to ``analytically continue'' these local expansions, i.e.,
how given a power of $p$ one can compute the rational functions which approximate
the entries in the Frobenius matrix to that power.

The essential point is that the theory guarantees that the holomorphic functions in the Frobenius
matrix $F(\Gamma)$ are ``overconvergent''. This implies that they converge on the $p$-adic projective line with
open disks of some unknown radius $s < 1$ removed around the poles of $r$. 
Let us notate this unknown larger domain by $D_s$. Assuming one
has an upper bound on $s$, and also an upper bound on the maximum
value $t$ taken by the Frobenius matrix on the closed set $D_s$, one can compute an upper bound
on the total degree of the rational functions needed to approximate $F(\Gamma)$ on this
domain modulo any given power of $p$. This upper bound allows one to determine how many
terms in the local expansion of $F(\Gamma)$ are required to compute the rational functions.
The knowledge of bounds $s$ and $t$ amounts to having {\it effective lower bounds}
on the $p$-adic decay of the entries in the matrix $F(\Gamma)$ (Definition \ref{Def-EPB}). 
We shall {\it assume} that
these effective lower bounds are {\it known}; for the explicit example we consider in Section 
\ref{Sec-AffineSurface}
we will explain exactly how to calculate them.

We refer the reader to \cite[Section 8.1]{LSH} for a detailed description of the relatively straightforward step of recovering the matrix of approximating rational functions from the local expansions given that these
bounds are known.

\section{An algorithm for computing $F:E^{1,n}_{2,rig} \rightarrow E^{1,n}_{2,rig}$}\label{Sec-MainAlg}

In this section we gather together the results from Sections \ref{Sec-RC}, \ref{Sec-RedAlgs}
and \ref{Sec-DefFrob} and present the main algorithm of the paper.

\subsection{Definitions and assumptions}\label{Sec-DefAss}

In this section we retain the definitions from Section \ref{Sec-CMD} and make
the assumption on the connection matrix from Section \ref{Sec-DefEetc}.
Specifically, we are given a pencil $X \rightarrow S$ of $k$-varieties such that:
\begin{itemize}
\item{The relative space $\E^\dagger:=H^n_{rig}(X/S) \cong H^n_{dR}(\gX_K/\gS_K) \otimes_A A^\dagger$ is a free $A^\dagger$-module of rank $m$.}
\item{The connection $\nabla^\dagger:\E^\dagger \rightarrow \E^\dagger \otimes \Omega^1_{A^\dagger}$ is given by a matrix $b(\Gamma)/r(\Gamma)$ with simple poles; here
$r(\Gamma) \in \Oh_K[\Gamma]$ is squarefree modulo $p$ and $d := \deg(r)$.}
\item{The space $E^{1,n}_{2,rig}$ is as in Definition \ref{Def-E2}.} 
\item{The morphism $\hat{\gX} \rightarrow \hat{\gS}$ is proper and smooth, so
that the base change theorem holds (Theorem \ref{Thm-Bert}).}
\item{The morphism $\gX_K \rightarrow \gS_K$ arises by extension of scalars from one
defined over an algebraic number field, so the local exponents are rational.}
\end{itemize}

\subsection{The comparison theorem in the geometric case}\label{Sec-CTholds}

\begin{theorem}\label{Thm-Cholds}
With definitions and assumptions as in Section \ref{Sec-DefAss}, Condition (O.C.) is met.
\end{theorem}

\begin{proof}
We have an overconvergent
$F$-isocrystal $(\E^\dagger,\nabla^\dagger,F)$
on $A^\dagger$. Choose any point $\Gamma = \gamma \in \Oh_{\bar{K}}$
such that $\ord_p(r(\gamma)) = 0$.
Dwork's trick of analytic continuation via Frobenius \cite[Prop. 3.1.2]{NKDwork}
tells us that the basis of location solutions to the differential
system $\nabla = 0$ converge on the open unit disk around $\Gamma = \gamma$.
The same is true for the second differential system $\check{\nabla} = 0$; consider
the ``dual $F$-isocrystal $(\E^\dagger,\check{\nabla}^{\dagger},F^{-1})$''
and apply Dwork's trick once again. So by Theorem \ref{Thm-Trans} condition
(O.C.) is met.
\end{proof}

\begin{theorem}\label{Thm-Bholds}
With definitions and assumptions as in Section \ref{Sec-DefAss}, condition (Rat.) is met.
\end{theorem}

\begin{proof}
Since $\gX_K \rightarrow \gS_K$  can be defined over the complex
numbers, it follows from the local monodromy theorem \cite[Thm (14.3)]{NK}.
\end{proof}

Theorems \ref{Thm-CT}, \ref{Thm-Cholds} and \ref{Thm-Bholds} together yield:

\begin{theorem}\label{Thm-CTholds}
With definitions and assumptions as in Section \ref{Sec-DefAss}., the comparison
theorem $E^{1,n}_{2,dR} \cong E^{1,n}_{2,rig}$ holds.
\end{theorem}

We note that our bounds on the growth of forms (Theorem \ref{Thm-Growth}) only holds when
the assumption (Rat.) is replaced by the stronger assumption (Prep.Rat.).
This will not always hold in the geometric case; however, we will give
examples in which it does hold (Section \ref{Sec-CompData}).

\subsection{Numerical approximations}

In this section we formalise the notion of a ``numerical approximation'' to a $p$-adic number.

We assume that elements in $\Oh_K$ are represented as $p$-adic expansions
with coefficients in some fixed set of representatives for the quotient
$\Oh_K/(p)$. Thus for any positive integer $N$, elements in the quotient
$\Oh_K/(p^N)$ can be represented in a unique manner via truncated $p$-adic
expansions. Elements in $\Oh_K[\Gamma,1/r(\Gamma)]/(p^N)$ have
an obvious representation via these truncated $p$-adic expansions.

\begin{definition}
Let $N$ be a positive integer, 
and $a \in A^\dagger =
K[\Gamma,1/r(\Gamma)]^\dagger$. Define
$N^\prime := N - \min(0,\ord_p(a))$.
A $p^N$-approximation to $a$ is a triple $(N,\ord_p(a),a_0)$ where
$a_0 \in \Oh_K[\Gamma,1/r(\Gamma)]/(p^{N^\prime})$ and
$a_0 - p^{-\min(\ord_p(a),0)}a = 0$ in $\Oh_K[\Gamma,1/r(\Gamma)]/(p^{N^\prime})$. 
\end{definition}

Thus taking $\hat{a}_0$ to be any preimage of $a_0$ in $\Oh_K[\Gamma,1/r(\Gamma)] \otimes
K$, we have that a $p^N$-approximation to $a$ defines an element
$a_1 := p^{\min(\ord_p(a),0)} \hat{a}_0$ such that $\ord_p(a - a_1) \geq N$. Conversely, 
given such an element $a_1$ one may canonically identify with it a $p^N$-approximation
to $a$. Intuitively a $p^N$-approximation amounts to knowledge ``modulo $p^N$''.

\subsection{Input/Output specification for the algorithm}\label{Sec-AlgInOut}

We retain the definitions and assumptions from Section \ref{Sec-DefAss}, and
now further assume that (Prep.Rat.) holds. So by Theorem
\ref{Thm-CTholds} our comparision theorem $E^{1,n}_{2,dR} \cong
E^{1,n}_{2,rig}$ holds, and moreover, we have effective bounds on the
growth of forms during reduction (Theorem \ref{Thm-Growth}).

 Let $N_I$ be a positive integer. 
We shall assume that we are given as {\it input} the following:
 \begin{itemize}
 \item{Input 1: The matrix $b(\Gamma)/r(\Gamma)$ for the connection $\nabla$.}
 \item{Input 2: A $p^{N_I}$-approximation to
 $F(\gamma)$ for one Teichm\"{u}ller specialisation $\Gamma = 
 \gamma$ of the matrix $F(\Gamma)$ for the action $F:H^n_{rig}(X/S) \rightarrow
 H^n_{rig}(X/S)$, i.e., an approximation to the $p$th power Frobenius action on
 $H^n_{rig}(X_{\bar{\gamma}})$ for some fibre $X_{\bar{\gamma}}$ of the family $X \rightarrow S$.} 
\end{itemize}
We also assume we are given (see Definition \ref{Def-EPB}):
\begin{itemize}
\item{Input 3: Effective $p$-adic bounds on the entries in the matrix $F(\Gamma)$.}
 \end{itemize}
 The algorithm gives as {\it output}:
 \begin{itemize}
 \item{Output: A $p^{N_O}$-approximation to a matrix for
$F:E^{1,n}_{2,rig} \rightarrow E^{1,n}_{2,rig}$,}
\end{itemize}
for some effectively computable $N_O < N_I$.
 The {\it loss of accuracy} is measured by
the difference $N_I - N_O$. We note in Section \ref{Sec-LossAcc} that
there exist effectively computable constants $\alpha^{\prime \prime},\beta^{\prime \prime}
 \geq 0$ such that one
may take $N_O := N_I - (\alpha^{\prime \prime} \log_p(N_I) - \beta^{\prime \prime})$, i.e., we have a ``logarithmic'' loss of accuracy. 

When $X$ is affine we have from equation (\ref{Eqn-Tsuzuki}) that 
$H^{n+1}_{rig}(X) \cong E^{1,n}_{2,rig}$ and
 so the matrix given as output yields an approximation to $p$th power Frobenius action on $H^{n+1}_{rig}(X)$.

\subsection{The algorithm}\label{Sec-TheAlg}

The algorithm comprises two steps

\begin{itemize}
\item{Step 1: From Inputs 1,2 and 3 use the ``deformation algorithm'' to compute a $p^N$-approximation to the matrix for $F:H^n_{rig}(X/S)
\rightarrow H^n_{rig}(X/S)$.}
\item{Step 2: From the output of Step 1 and Input 1, use the algorithm from Section
\ref{Sec-RedAlgs}
 to compute a $p^{N_O}$-approximation matrix for $F:E^{1,n}_{2,rig} \rightarrow E^{1,n}_{2,rig}$.}
\end{itemize}
The intermediate precision $N$, where $N_O < N < N_I$, 
can be computed from the input data, see Section \ref{Sec-LossAcc}.

\subsubsection{Step 1}

This step was described in detail in Sections \ref{Sec-LocalDef} and \ref{Sec-GlobalDef}. (One may
need to make a change of basis so that $\gamma = 0$.)
Estimates from Input 3 determine the $\Gamma$-adic accuracy required in the local
deformation to compute an $p^N$-approximation to the
global matrix $F(\Gamma)$ itself, whose entries lie in $A^\dagger$.

\subsubsection{Step 2}

Let $\B$ be a set of forms in $A^m \otimes d\Gamma$ whose image in
$E^{1,n}_{2,rig}$ are a $K$-basis for this space. This may be determined
by computing the set of exponents $E(\nabla,\B)$, then performing the linear algebra computation described in the
proof of Theorem \ref{Thm-BasisdR}. Theorem \ref{Thm-CT} assures us that
this set gives a basis for $E_{2,rig}^{1,n}$. We can assume that $\ord_p(e) = 0$
for all $e \in \B$.

For each $e \in \B$, one computes a $p^N$-approximation of
the image \linebreak $F(\Gamma)p\Gamma^{p-1} e^\sigma(\Gamma^p)d\Gamma$. Then one performs a radix
conversion from $\Gamma$-adic to $r(\Gamma)$-adic expansions, so that the
input is in the appropriate form for the reduction algorithm. (It is actually much better
in practice to compute an $r$-adic expansion of the matrix
$F(\Gamma)/r^\sigma(\Gamma^p)$, and recover the $r$-adic expansions for each
image form via a $r$-adic multiplication routine.
 It turns out that these radix conversions are very time consuming, so one wishes to minimize the number performed.)
Then use the reduction
algorithm from the proof of Theorem \ref{Thm-FGEdR}, plus the final linear
algebra step from the proof of Theorem \ref{Thm-BasisdR}, to write
this as a $p^{N_O}$-approximation to a $K$-linear combination of elements in $\B$ plus an
$p^{N_O}$-approximation to an element in $\nabla^\dagger(\E^\dagger)$.

\subsection{Analysis}

\subsection{Loss of accuracy}\label{Sec-LossAcc}

Theorem \ref{Thm-Growth} and the analysis following Theorem \ref{Thm-DiffSysLoss} together show that the matrix computed
in the algorithm is a $p^{N_O}$-approximation to a matrix for the action of
$F:E^{1,n}_{2,rig} \rightarrow E^{1,n}_{2,rig}$ for $N_O$ such that
\[ N_O = N_I - (\alpha^{\prime \prime} \log_p(N_I) + \beta^{\prime \prime}).\]
Here $\alpha^{\prime \prime},
\beta^{\prime \prime} \geq 0$ are constants which may be computed from $m$, $p$, the local
exponents of the connection, and the effective $p$-adic bounds on $F(\Gamma)$ (Input 3).
We shall not present explicit formulae for $\alpha^{\prime \prime}$ and 
$\beta^{\prime \prime}$ in the general
case as they are rather complicated. We note that the discussion following
Theorem \ref{Thm-DiffSysLoss} allows
one to compute the intermediate precision $p^N$ which is attained after
Step 1.

\subsection{Time and space complexity}

The time and space complexity may be calculate given the effective $p$-adic bounds on the
matrix $F(\Gamma)$, and also a bound on the height of the local monodromy eigenvalues, c.f.
Section \ref{Sec-RT}. 
 We do not present an explicit expression for the general case
since it is rather complicated. Let us just make a few observations on Step 1:
The calculation of the local solution matrix in Step 1 (Section \ref{Sec-LocalDef})
is fast, both in theory and practice, since it just requires the iteration of a short 
linear recurrence; however, Method 1 is rather space consuming in comparison to 
Method 2. The analytic continuation step requires only a single multiplication
by a power of $r(\Gamma)$ (computed modulo a power of $\Gamma$). The
radix conversion, though in theory ``quasi-linear time'' \cite[Alg. 9.14]{vGG}, is in practice
rather time consuming.

\section{The Frobenius matrix of an affine surface}\label{Sec-AffineSurface}

In this section we apply the algorithm in Section \ref{Sec-MainAlg} to compute
to any required numerical precision a matrix for the $p$th power Frobenius map
acting on the middle dimensional rigid cohomology of a certain affine surface. Specifically,
we consider an open subset $X$ of the affine surface defined by an equation of the
form $Z^2 = \bar{Q}(X,\Gamma)$, subject to certain smoothness assumptions. The algorithm
from Section \ref{Sec-TheAlg} allows the efficient computation of an approximation
to the Frobenius map $F: H^2_{rig}(X) \rightarrow H^2_{rig}(X)$, provided one can
obtain the auxiliary inputs 1, 2 and 3 (Section \ref{Sec-AlgInOut}). After defining the
surface in Section \ref{Sec-ExplicitPencil}, we describe how the necessary
auxiliary inputs may be calculated (Section \ref{Sec-AuxInput}). Having specified
some local monodromy restrictions to ensure applicability of the algorithm in 
Section \ref{Sec-MainAlg} (see Section \ref{Sec-LMA}), we then give a precise complexity
analysis (Theorem \ref{Thm-TimeSpace}). 

In Section \ref{Sec-CompactSurface} we shall apply the results of the present section to
compute the full zeta function of a compactification $\bar{X}$ of the open surface $X$. We 
report on our Magma implementation of this final algorithm in Section \ref{Sec-CompData}.

We retain the notation in Section \ref{Sec-RC}.  In particular, recall that $k = \fq$ is the finite field with
$q$ elements, and $K$ the unramified extension of $\Qp$ of degree $[k:\fp]$. We assume now
that the characteristic $p$ is odd.
The ring
of integers of $K$ is denoted $\Oh_K$.
Let us further assume we are given $L \supseteq \Q$ an algebraic number field
with ring of integers $\Oh_L$. Assume that we have an embedding
$\Oh_L \subset \Oh_K$ and $\Oh_L/(p) \cong \Oh_K/(p) \cong \fq$. 

\subsection{Definition of the pencil}\label{Sec-ExplicitPencil}

Let $Q(X,\Gamma) \in \Oh_L[X,\Gamma]$ and denote by $\bar{Q}(X,\Gamma) \in k[X,\Gamma]$
its reduction modulo $p$. We shall assume that both $Q$ and $\bar{Q}$ are monic in $X$ of
degree $2g+1$ where $\gcd(p,2g+1) = 1$.
Let $\res(\Gamma) := \Res(X,Q,\frac{\partial Q}{\partial X}) \in \Oh_L[\Gamma]$ be the
Sylvester resultant w.r.t. $X$ of $Q$ and $\frac{\partial Q}{\partial X}$, see \cite[Pages 150-151]{CLO} or
Section \ref{Sec-ComputeNabla}. (The notation $r(\Gamma)$ is reserved for
the monic denominator of the connection matrix, which is a factor of $\res$, see
Section \ref{Sec-ComputeNabla}.) Assume that $\res(\Gamma)$ has leading coefficient
a unit modulo $p$, and $\res(0) \ne 0 \bmod{p}$; in particular, it does not vanish identically modulo $p$.
Define the $\Oh_K$-schemes
\[ \gX:=\Spec(B) \mbox{ where } B:= \Oh_K[X,\Gamma,Z,\res(\Gamma)^{-1}]/(Z^2 - Q(X,\Gamma))\]
and
\[ \gS:=\Spec(A) \mbox{ where } A:=\Oh_K[\Gamma,1/\res(\Gamma)].\]
Let $\bar{B}:= B \otimes_{\Oh_K} k$ and $\bar{A} := A \otimes_{\Oh_K} k$ be the reduction of the
coordinate rings modulo $p$. Define the $k$-schemes
\[ X := \Spec(\bar{B}) \mbox{ and } S:= \Spec(\bar{A}).\]
We have the obvious commutative diagrams
\[
\begin{array}{ccccccc}
A & \hookrightarrow & B & & \gX & \rightarrow & \gS\\
\downarrow & & \downarrow & \mbox{ and } & \downarrow & & \downarrow\\
\bar{A} & \hookrightarrow & \bar{B} & & X & \rightarrow & S
\end{array}
\]
where the vertical maps in the second diagram take special fibres.  
Recall that the generic fibres are denoted $\gX_K$ and $\gS_K$, respectively.
The horizontal maps in the second diagram are smooth morphisms of smooth schemes,
and the fibres are (affine) hyperelliptic curves.

The relative cohomology spaces which concern us are:
\[ \E := H_{dR}^1(\gX_K/\gS_K) = \left\langle \frac{X^i dX}{\sqrt{Q}}\,|\,0 \leq i < 2g \right\rangle_{A}\]
and
\[ \E^\dagger := H_{rig}^1(X/S) = \left\langle \frac{X^i dX}{\sqrt{Q}}\,|\,0 \leq i < 2g 
\right\rangle_{A^\dagger}.\]
Here $\sqrt{Q}$ denotes the image of $Z$ in $B$ (precisely, in $B^\dagger$ for the second
space). We refer the reader to Section \ref{Sec-Pencils} for a description of the ring
$A^\dagger$. 
That $H_{dR}^1(\gX/\gS)$ and $H_{rig}^1(X/S)$ are spanned by these
forms is shown in \cite[Secs. 4.2, 5.4]{LFFA}. That they form a basis follows by a specialisation
argument, and the fact that the dimension of the first de Rham (rigid resp.) cohomology space of
any fibre in the family $\gX \rightarrow \gS$ ($X \rightarrow S$ resp.) is $2g$. Alternatively, 
see \cite{HH}. Note that we do not need to appeal to the finiteness and comparison
theorems in Section \ref{Sec-RCT}, since we can establish the necessary results directly.

\subsection{The spectral sequence}

The next proposition shows that the algorithm in Section \ref{Sec-MainAlg} in the present case
computes a numerical approximation to a matrix for $F:H^2_{rig}(X) \rightarrow H^2_{rig}(X)$.

\begin{proposition}\label{Prop-SSdeg}
With the morphism $X \rightarrow S$ as defined in Section \ref{Sec-ExplicitPencil}, we have
\[ H^2_{rig}(X) \cong E^{1,1}_{2,rig}.\]
\end{proposition}

\begin{proof}
Follows from equation (\ref{Eqn-Tsuzuki}) since $X$ is affine and $\Omega^1_{R/\V}$ is free
of rank one; here $R:=\Gamma(\gS,\Oh_\gS)$ and $\V := \Oh_K$.
\end{proof}

\subsection{Auxiliary data: Inputs 1,2 and 3}\label{Sec-AuxInput}

We now explain how to compute the auxiliary information needed as input to our main 
algorithm (Section \ref{Sec-MainAlg}) in the case of the surfaces presently under consideration.

\subsubsection{Input 1: the matrix for $\nabla$}\label{Sec-ComputeNabla}

The connection $\nabla$ acts on $\E$ by taking the derivative w.r.t. $\Gamma$ of the basis
elements $X^i dX/\sqrt{Q}$, and then applying the reduction algorithm of Kedlaya
to write the image as a linear combination over $A$ of the basis elements plus
$\frac{d}{dX}(g) \otimes d\Gamma$ for some element $g \in \E$. We now give
explicit formulae for computing the matrix for the connection (based upon
Magma code written by the author).

For an element $a \in L(\Gamma)[X]$ and $i \in \Z$, denote by $\Coeff(a,i)$ the coefficient
of $X^i$ in $a$. Let $\delta := 2g+1$. Let $M$ be the Sylvester matrix w.r.t. X of $Q$ and $\frac{\partial Q}
{\partial X}$. Explicitly, $M \in M_{2\delta-1}(L[\Gamma])$ and for $1\leq j \leq 2\delta-1$
\[ M_{ij} := 
\left\{
\begin{array}{ll}
\Coeff(X^{\delta-1-i}Q,2\delta-1-j) & \mbox{ for $1 \leq i \leq \delta-1$}\\
\Coeff(X^{2\delta-1-i}\frac{\partial Q}{\partial X},2\delta-1-j) & \mbox{ for $\delta \leq i \leq 2\delta-1$}.
\end{array}
\right.
\]
We have assumed the determinant of this matrix $\res(\Gamma)$ (Sylvester resultant) is non-zero
modulo $p$. Define $E$ to be the $\delta-1 \times 2\delta-1$ matrix over $L[\Gamma]$ with
\[ E_{ij} := -\frac{1}{2} \Coeff\left(X^{i-1}\frac{\partial Q}{\partial \Gamma},(2\delta-1)-j\right).\]
Let $F := EM^{-1}$, a $\delta-1 \times 2\delta-1$ matrix over $L[\Gamma,1/\res(\Gamma)]$. Let the vectors
$a,b,c \in L[\Gamma,1/\res(\Gamma)][X]^{\delta-1}$ be defined as follows: For $1 \leq i \leq \delta-1$
\[
\begin{array}{rcl}
a_i & := & \sum_{j = 1}^{\delta-1} F_{i,\delta-j} X^{j-1}\\
b_i & := & \sum_{j = 1}^{\delta-1} F_{i,2\delta-j} X^{j-1}\\
m_i & := & a_i + 2\frac{\partial b_i}{\partial X}.
\end{array}
\]
Then the connection matrix $B(\Gamma) \in M_{\delta-1}(L[\Gamma,1/\res(\Gamma)])$ is defined by
$B_{i,j} := \Coeff(m_j,i-1)$. One can uniquely write $B(\Gamma) = b(\Gamma)/r(\Gamma)$ where
$r|\res$, and $r$ is monic and coprime to some entry in the matrix $b \in M_{\delta - 1}(L[\Gamma])$.

We shall impose some restrictions on the connection matrix $B(\Gamma)$ in Section \ref{Sec-LMA}.

 \subsubsection{Input 2: The Frobenius matrix of a fibre}
 
We take the fibre at $\gamma = 0$, noting $\res(0) \ne 0 \bmod{p}$.
The Frobenius matrix of the fibre $Z^2 = \bar{Q}(X,0)$ can be computed
using Kedlaya's original algorithm \cite{KK1}; the implementation by Michael Harrison is available
with the documentation accompanying the Magma program.

\subsubsection{Input 3: Effective $p$-adic bounds for $F:\E^\dagger \rightarrow \E^\dagger$}\label{Sec-EffBounds}
 
The Frobenius matrix $F(\Gamma)$ can in principle be calculated by applying Kedlaya's
algorithm to the ``generic'' hyperelliptic curve in the family, which is defined over the function field $\fq(\Gamma)$. From the point of view of complexity theory this is not a good idea; it is faster to use the indirect method of the ``deformation algorithm''. However, this direct method is a good way to calculate
effective $p$-adic bounds for the matrix $F(\Gamma)$.

Specifically, fix $i,j$ with $1 \leq i,j \leq 2g$. Let $f(\Gamma)$ be the $(i,j)$th entry
in the matrix $F(\Gamma)$. Then $f(\Gamma)$ is the coefficient of $X^i dX/\sqrt{Q}$ in the
expression one obtains by reducing the form
\[ \sigma\left(\frac{X^j dX}{\sqrt{Q}}\right) =
\frac{pX^{p(j+1)-1}}{Q^{p/2}}\left(1 - \frac{Q^p - Q^\sigma(X^p,\Gamma^p)}{Q^p} \right)^{-1/2}dX\]
using the ``generic'' version of Kedlaya's algorithm.
Here $\sigma$ is the map sending $\Gamma \mapsto \Gamma^p$, $X \mapsto X^p$,
and acting like the $p$th power Frobenius automorphism on $K$.
 We can write
$Q^p - Q^\sigma(X^p,\Gamma^p) = p R(X,\Gamma)$ for some unique $R \in \Oh_K[X,\Gamma]$ with
\[ \deg_X(R) < p\deg_X(Q),\,\deg_\Gamma(R) \leq p\deg_\Gamma(Q)\]
with second inequality strict if $Q$ is monic in $\Gamma$.
Then
\[ \sigma(X^j dX/\sqrt{Q}) = \sum_{\ell  = 0}^\infty {-1/2 \choose \ell} p^{\ell+1} \frac{X^{p(j+1) - 1}R^\ell}{Q^{p(\ell+(1/2))}}dX.\]
The $\ell$th term in this series can be reduced modulo exact forms using $p\ell + \lfloor \frac{p}{2}
\rfloor$ applications
of Kedlaya's ``pole reduction formula'', see \cite[Section 4.2]{LFFA}. Each application requires one division by the resultant
$\res(\Gamma)$. By an easy specialisation argument, \cite[Lemma 2]{KK1} implies that reduction
of the $\ell$th term requires a cumulative division by at most $p^{\lfloor \log_p(p(2\ell + 1)) \rfloor}$.

Write
\[ f(\Gamma) = \sum_{k = -\infty}^{\infty} f_k(\Gamma) \res(\Gamma)^i\]
where $f_k \in K[\Gamma]$ with $\deg(f_k) < \deg(\res)$. The argument in the preceding paragraph implies the following.

\begin{proposition}\label{Prop-Eff1}
For $k < 0$ we have the lower bound
\[ \ord_p(f_k(\Gamma)) \geq  (\ell + 1) - \lfloor \log_p(p(2\ell+1)) \rfloor \]
where $\ell$ is the smallest integer such that $p\ell + \lfloor \frac{p}{2} \rfloor \geq |k|$.
(Explicitly, $\ell:= \lfloor \frac{2|k| - p + 1}{2p} \rfloor$ and so $\ord_p(f_k(\Gamma)) \geq 
\lfloor |k|/p \rfloor - \lfloor \log_p(2|k| + 1) \rfloor$.)
\end{proposition}
 
A lower bound for $k \geq 0$ requires a more detailed analysis: Let $\Adj(M)$ be the adjoint of the
Sylvester matrix. Each application of Kedlaya's pole reduction formula increases
the degree in $\Gamma$ (degree of numerator minus degree in denominator)
by $\deg_\Gamma(\Adj(M)) - \deg_\Gamma(\res)$. The degree in $\Gamma$ of the numerator $X^{p(j+1)-1}R^\ell$
 of the $\ell$th term in the series is $\ell \deg_\Gamma(R) < \ell (p\deg_\Gamma(Q) - 1)$.
Thus after $p \ell + \frac{p-1}{2}$ applications of Kedlaya's pole reduction formula, the degree
in $\Gamma$ of the reduction of the $\ell$th term in the series is at most
\[ \kappa(\ell) := \left(p\ell + \frac{p-1}{2}\right)\left(\deg_\Gamma(\Adj(M)) - \deg_\Gamma(\res)\right) + \ell \deg_\Gamma(R).\]
We note that  the modest use of Kedlaya's formula for reducing the ``pole at infinity'' required in the
calculation of $F(\Gamma)$ does not increase the degree in $\Gamma$ (or introduce
powers of $\tilde{r}(\Gamma)$ on the denominator).
Thus we deduce:

\begin{proposition}\label{Prop-Eff2}
For $k \geq 0$ we have the lower bound
\[ \ord_p(f_k(\Gamma)) \geq (\ell + 1) - \lfloor \log_p(p(2\ell + 1)) \rfloor\]
where $\ell$ is the smallest integer such that $\frac{\kappa(\ell)}{\deg_\Gamma(\res)} \geq k$. (If
no such $\ell$ exists then the term $f_k(\Gamma)$ is zero.) 
\end{proposition}

Explicitly, define 
\begin{equation}\label{Eqn-Delta}
 \delta:= \frac{\deg_\Gamma(\Adj(M))}{\deg_\Gamma(\res)},\,
\delta^\prime :=\frac{\deg_\Gamma(R)}{p\deg_\Gamma(\res)} \leq \frac{\deg_\Gamma(Q)}
{\deg_\Gamma(\res)}.
\end{equation}
Then assuming $\delta + \delta^\prime
\geq 1$ one takes $\ell$ the floor of $\frac{2k - (p-1)(\delta - 1)}{2p (\delta + \delta^\prime - 1)}$.

We now state a conjecture to which we shall refer later.

\begin{conjecture}\label{Conj-FinitePoles}
The Frobenius matrix $F(\Gamma)$ has a pole of finite order at infinity, rather than
an essential singularity.
\end{conjecture}

Conjecture \ref{Conj-FinitePoles} thus claims that $f_k(\Gamma) = 0$ for sufficiently
large positive $k$.

\subsection{Local monodromy assumptions}\label{Sec-LMA}

In this section we state some further restrictions made to ensure that the conditions
required for the application of the main algorithm in Section \ref{Sec-MainAlg} are met.
Specifically, we need that the connection matrix $B(\Gamma)$ from Section
\ref{Sec-ComputeNabla} is of the form required
in Section \ref{Sec-DefAss}, and that condition (Prep.Rat.) is met. To simplify the complexity analysis and
to keep in line with our actual implementation in Section \ref{Sec-CompData} we shall
in fact make stronger assumptions, as follows.

Recall that $Q \in \Oh_L[X,\Gamma]$ with 
$2g + 1:=\deg_X(Q)$, that $\tilde{r}(\Gamma)$ is the
Sylvester resultant of $Q$ and $\frac{\partial Q}{\partial X}$ w.r.t. $X$, and $r$ the monic
factor of $\tilde{r}$ which is the denominator of the connection matrix $B(\Gamma) = 
b(\Gamma)/r(\Gamma)$ when in lowest terms. Define $h:=\deg_\Gamma(Q)$.

We {\it assume} that $r(\Gamma) \bmod{p}$ is squarefree, and that the Laurent expansion of $B(\Gamma)$
has only negative terms. We say then that $B(\Gamma)$ has {\it only simple poles modulo $p$}.
This ensures that the algorithm in the proof of Theorem
\ref{Thm-FGEdR} works. Let us {\it assume} that the local monodromy eigenvalues around each
singular point are prepared, so condition (Prep.Rat.) is met and we may apply the precision
loss bounds in Theorem \ref{Thm-Growth}. 

To obtain a nice basis for $E^{1,1}_{2,rig}$,
let us further {\it assume} that
the local monodromy around the finite poles is nilpotent, and that zero does not
occur as a local monodromy eigenvalue around the pole at infinity. 
 In this case we
may take as our basis for $E^{1,1}_{2,rig}$ the elements $\{b_{ik}\}$ where
$0 \leq i \leq d-2$ $(d :=\deg(r))$ and $1 \leq k \leq  2g$, and the element $b_{ik} \in H^1_{rig}(X/S)$ is the column vector with zeros in positions $j \ne k$, and in position $j = k$ the
$1$-form $\Gamma^i d\Gamma/r(\Gamma)$. Note that the dimension of this
space is $2g (d-1)$.

For the complexity analysis, we shall need bounds on the height of the local monodromy eigenvalues.
Let us {\it assume} that a common denominator for the local monodromy eigenvalues around
each singular point is $2(2g+1)$, and when written w.r.t. this denominator the numerator does
not exceed $h(2g-1)$ in absolute value. Under the assumption that $B(\Gamma)$ has only
simple poles modulo $p$, we believe that one may prove that the bound on the denominator
should always holds by a topological argument. Likewise, the author expects that the bound
on the numerator should also hold, although offers no proof of this. 

\begin{note}\label{Note-LME}
We point out that ``generically'' in any nice family of polynomials both
$r(\Gamma)$ and $\tilde{r}(\Gamma)$ are squarefree and have equal degree. 
In this case,
one observes experimentally, and expects to be able to prove, that all residue matrices around finite poles are nilpotent. However, the assumption that the degree in $\Gamma$
of the connection matrix $B(\Gamma)$ is less than zero does {\it not} hold generically. For any family of
polynomials $Q$, e.g. with fixed Newton polytope, one can calculate restrictions on the
coefficients which must be met. The author has no idea of the geometric significance
of this assumption. When the assumption does hold, the local exponents at infinity
are observed to exhaust the set $\left\{\frac{\pm j h}{2(2g+1)}| 1 \leq j \leq 2g-1,\, j \mbox{ odd}\right\}$.
\end{note}

\subsection{Analysis}

We shall use soft-Oh notation, to hide logarithmic factors in the time and space complexity
\cite[Def. 25.8]{vGG}.

Let $N_O$ be  a positive integer which depends upon the equation
$Z^2 = \bar{Q}(X,\Gamma)$ in some manner  --- we shall specify precisely
how later. Assume that 
\begin{equation}\label{Eqn-GrowthNO}
g^2\left(1 + \frac{h}{p} + \log_p(gh)^2 \right) = \Oh(N_O),
\end{equation}
i.e., the integer $N_O$ grows at least as fast as the expression on the
lefthand-side as $g,h$ and $p$ vary.

\subsubsection{Numerical approximations}\label{Sec-NumApp}

Assume that one wishes to compute a $p^{N_O}$-approximation to the
$p$th power Frobenius matrix $F: H^2_{rig}(X) \rightarrow H^2_{rig}(X)$.
Then Theorem \ref{Thm-Growth}, equation (\ref{Eqn-AbsEst}) in
Note \ref{Note-Constants}, inequality (\ref{Eqn-M1loss}), Propositions \ref{Prop-Eff1} and \ref{Prop-Eff2}, and the local monodromy assumptions in Section \ref{Sec-LMA} show that it suffices to take the
initial $p$-adic accuracy $N_I$ such that $N_I  - (\alpha^{\prime \prime}
 \log_p(N_I) + \beta^{\prime \prime}) \geq N_O$
for some effective constants $\alpha^{\prime \prime},
\beta^{\prime \prime} \geq 0$. For implementations one needs
to compute the loss of accuracy precisely; however for our complexity estimates it
is enough to observe
$\alpha^{\prime \prime} = \Oh( g^2 \log_p(g))$ and 
$\beta^{\prime \prime} = \Oh( g^2 (1 + (h/p) + \log_p(gh)^2))$.
Here are more details. First, Propositions \ref{Prop-Eff1} and \ref{Prop-Eff2} combined
with the observation $(\delta + \delta^\prime - 1)\deg_\Gamma(\res) = \Oh(gh)$ shows
the following: the $\Gamma$-adic accuracy needed in solution of the differential system
in Step 1 is $\Oh(pgh N_I)$. From 
inequality (\ref{Eqn-M1loss}), the loss of accuracy in this step is $\Oh(g\log_p(pgh N_I)) = 
\Oh(g (\log_p(N_I) + \log_p(gh)))$. Second, the maximum pole order encountered in Step 2
 is $\Oh(p ghN)$ where $N < N_I$ is the intermediate accuracy, so the loss of accuracy
 in Step 2 is $\Oh(\alpha \log_p(pghN) + \beta) = \Oh (\alpha \log_p(N_I) + 
 \alpha\log_p(gh) + \beta)$ where $\alpha,\beta$ are the constants in
 Theorem \ref{Thm-Growth}. From equation (\ref{Eqn-AbsEst}) we
 have $\alpha = \Oh( g^2 \log_p(g))$ and $\beta = \Oh( g^2 (1 + (h/p) + \log_p(gh))$. 
 Our claim on the loss of accuracy now follows. Moreover, from equation
 (\ref{Eqn-GrowthNO}) we see that the initial $p$-adic accuracy $N_I$ satisfies
 $N_I = \tilde{\Oh}(N_O)$.

\subsubsection{Time and space complexity}\label{Sec-RT}

We now give a precise complexity analysis of the time and space required to 
compute a numerical approximation to the $p$th power Frobenius
matrix $F:H^2_{rig}(X) \rightarrow H^2_{rig}(X)$ using the algorithm
in Section \ref{Sec-MainAlg}.

\begin{theorem}\label{Thm-TimeSpace}
Let the affine surface be defined as in Section \ref{Sec-ExplicitPencil}, and 
assume the local monodromy conditions specified in Section \ref{Sec-LMA} hold.
We recall that $X$ is an open subset of the smooth surface defined
by the equation $Z^2 = \bar{Q}(X,\Gamma)$ over the field $\fq$ of characteristic
$p$, and $2g + 1:=\deg_X(Q),\,h :=\deg_\Gamma(Q)$.
Let the positive integer $N_O$ satisfy the growth condition (\ref{Eqn-GrowthNO}).
Then one may compute a $p^{N_O}$-approximation to the $p$th power
Frobenius matrix $F: H^2_{rig}(X) \rightarrow H^2_{rig}(X)$ via the
algorithm in Section \ref{Sec-MainAlg} in $\tilde{\Oh}(N_O^2 g^5 h^2 p \log(q))$ bit 
operations, using $\tilde{\Oh}(N_O^2 g^3 h p\log(q))$ bits of space.
\end{theorem}

\begin{proof}
Since $N_O$ satisfies growth condition (\ref{Eqn-GrowthNO}), from Section
\ref{Sec-NumApp} we see that the initial $p$-adic accuracy $N_I$ satisfies
$N_I = \tilde{\Oh}(N_O)$. For the purposes of the complexity analysis, we shall
forget about the intermediate accuracy $N$, with $N_O < N < N_I$ mentioned in 
Section \ref{Sec-TheAlg}, and just assume we work with $p^{N_I}$-approximations
throughout the algorithm.

Step 1: Using the estimates from Propositions \ref{Prop-Eff1} and \ref{Prop-Eff2},
we see that the $\Gamma$-adic accuracy required in Step $1$ is
$\Oh(N_I p \mu)$, where 
\[ \mu :=\max\{\deg_\Gamma(Q),\deg_\Gamma(r),\deg_\Gamma(\Adj(M))\} = \Oh(hg).\]
We consider the time/space required to compute an approximation to
$C(\Gamma)$ in Section \ref{Sec-Method1}:
The coefficients of $\Gamma$ are $g \times g$ matrices, whose entries are
$p^{N_O}$-approximations of elements of the $p$-adic field $K$. Moreover, the
growth bounds given in the analysis following Theorem \ref{Thm-DiffSysLoss} shows that
each coefficient requires $\tilde{\Oh}(\log(q)N_I)$ bits of space. This gives a space
requirement of $\tilde{\Oh}(N_I^2 g^3 h p \log(q))$ bits. For the time, we observe that
recurrence (\ref{Eqn-ApproxSol}) has length bound by $\max\{\deg(b)+1,\deg(r)\} = 
\Oh(gh)$, and involves multiplication of $g \times g$ matrices. Thus the time
to compute an approximation to $C(\Gamma)$ 
is $\tilde{\Oh}(N_I^2 g^{2 + \omega} h^2 p \log(q))$ bit operations.
One may further compute the approximation to the local Frobenius matrix $F(\Gamma)$
in this time/space, using (\ref{Eqn-F}).
Using the fast radix conversion algorithm in \cite[Alg 9.14]{vGG}, these time and
space estimates are enough for the
the analytic continuation and radix conversion
steps required to make the input suitable for Step 2.

Step 2: The matrix $F$ has size $(d-2)2g = \Oh(g^2 h)$, where $d = \deg(r(\Gamma))$.
Thus $\Oh(g^2 h)$ applications of the reduction algorithm from Section \ref{Sec-RedAlgs}
are required.
It is time saving in terms of the parameter $g$ to precompute the
inverses  ``$(-\ell r^\prime I_m + b)^{-1}, (\ell^\prime I_m + b_{d-1})^{-1}$''  for
$\ell$ and $\ell^\prime$ in the necessary ranges, as these do not depend 
on the element being reduced. The number of the former inverses
is $\Oh(N_I  p)$ and each inverse takes $\tilde{\Oh}(N_I \log(q) \times gh \times g^3)$ bit operations
to compute; the factor $gh$ arising since the inverse is computed modulo
the polynomial $r(\Gamma)$ which has degree $\Oh(gh)$. There are $\Oh(N_I p gh)$ of the
latter inverses to compute, but each only requires $\Oh(N_I \log(q) g^3)$ bit operations.
Thus precomputation of the matrix inverses takes
$\tilde{\Oh}(N_I^2 g^4 h p \log(q))$ bit operations and one needs
$\Oh(N_I^2 g^3 h p \log(q))$ bits to store them. The reduction of finite poles
requires $\Oh(N_I p)$ steps, and each step taking $\tilde{\Oh}(N_I g^3 h \log(q))$ bit operations;
reduction of the pole at infinity requires $\Oh(N_I p gh)$ steps, but each step only
taking $\tilde{\Oh}( N_I g^2 \log(q))$ bit operations. The time for the reduction of forms
is thus $\Oh(N_I^2 g^5 h^2 p \log(q))$ bit operations, and this step requires
$\tilde{Oh}(N_I^2 g^3 p \log(q))$ bits of space. (The time without precomputation of matrix
inverses would be $\tilde{\Oh}(N_I^2 g^6 h^2 p \log(q))$ bit operations.) This completes
the proof.
\end{proof}

\section{The zeta function of a compact surface}\label{Sec-CompactSurface}

This section is a direct continuation of Section \ref{Sec-AffineSurface}. In particular, 
throughout this section we retain the definitions and assumptions
in the preamble to that section, as well as those given in Sections \ref{Sec-ExplicitPencil} and
\ref{Sec-LMA}.

\subsection{The zeta function of the open surface}\label{Sec-P1}

In this section we consider the zeta function $Z(X,T)$ of the smooth affine surface $X$
over $\fq$.
The trace formula in rigid cohomology for smooth affine varieties shows
\begin{equation}\label{Eqn-ZetaX}
 Z(X,T) = \frac{P_1(X,T)}{P_2(X,T) P_0(X,T)}
 \end{equation}
where $P_i(X,T) :=\det(1 - T q^2 F^{-\log_p(q)}|H^i_{rig}(X)) \in 1 + T\Qp[T]$.
Certainly $H^0_{rig}(X)$ is a one-dimensional $\Qp$-vector space, and
$P_0(X,T) = (1 - q^2 T)$.

\begin{proposition}\label{Prop-P1}
Let the polynomial $P_1(S,T)$ be the numerator of the zeta function of the
open subset $S$ of the projective line; so $P_1(S,T)$ is a product of
cyclotomic polynomials. Then $P_1(X,T) = P_1(S,qT) \in 1 + T \Z[T]$.
\end{proposition}

\begin{proof}
It is enough to consider the terms $E_{2,rig}^{0,1}$ and $E_{2,rig}^{1,0}$ in the spectral
sequence for $X \rightarrow S$, c.f. Section \ref{Sec-Tsuzuki} and
\cite[Eqn (17)]{KO}.
We have
\[ E_{2,rig}^{0,1} := \ker(\nabla^\dagger) \cong \ker(\nabla).\]
The isomorphism follows from \cite[Cor. 2.6]{BC}. We claim the latter space is
zero-dimensional: Let $v \in \ker(\nabla)$. Recalling from Section \ref{Sec-LMA} that the local monodromy eigenvalues around finite poles are all zero, 
expanding $v$ around the finite poles
one deduces that $v \in K^{2g}$. Since the local monodromy eigenvalues around the pole at infinity
are non-zero, expanding $v$ around this pole one deduces $v = 0$. We have
\[ E_{2,rig}^{1,0} := \coker\left( \frac{d}{d\Gamma} : H^0_{rig}(X/S) \rightarrow H^0_{rig}(X/S)d\Gamma
\right).\]
But $H^0_{rig}(X/S) \cong A^\dagger$, the weak completion of the coordinate ring of
$S$. So 
\[ \det(1 - Tq^2 F^{-\log_p(q)}| E_{2,rig}^{1,0}) =  \det(1 - (Tq) qF^{-\log_p(q)}| H_{rig}^1(S)) = 
P_1(S,qT).\]
\end{proof}

\begin{proposition}
The polynomial $P_2(X,T)$ has integer coefficients.
\end{proposition}

\begin{proof}
Integrality follows from Proposition \ref{Prop-P1} since the zeta function itself is
a power series with integer coefficients.\end{proof}

Kedlaya's $p$-adic analogue of Deligne's main theorem tells us that the complex
absolute values of reciprocal zeros of $P_2(X,T)$ belong to the set $\{1,q^{1/2},q\}$
\cite{KKII}; we will deduce this is an elementary manner in Section \ref{Sec-Compact}.

\subsection{The zeta function of a compactification}\label{Sec-Compact}

In this section we show that the full zeta function $Z(\bar{X},T)$ of a compactification $\bar{X}$ of
$X$ may be easily recovered given the first $\Oh(gh)$ coefficients
of $P_2(X,T)$ to precision modulo $p^{N}$ where $N = \Oh(gh \log(q))$.

To simplify the analysis, and keep in line with our actual implementation, we shall
make some further restrictions on the polynomial $Q(X,\Gamma)$. Recall that
$Q$ is monic in $X$ of degree $2g+1$, and has degree $h$ in $\Gamma$. Let
us further assume that it is monic in $\Gamma$ with $h$ odd, has constant term
$1$, and that $2g+1,\,h$ and the prime $p$ are mutually coprime.
Moreover, assume that all other terms in $Z^2 - Q(X,\Gamma)$ have exponents
lying within or on the boundary of the polytope $\Delta$ with vertices the origin and the points
$(2g+1,0,0),(0,h,0)$ and $(0,0,2)$. 
Then the Newton polytope of $Z^2 - Q(X,\Gamma)$ (taken modulo any prime number
$p$) is the simplex $\Delta$. We assume that $Z^2 - \bar{Q}(X,\Gamma)$ is {\it non-degenerate}
w.r.t. the polytope $\Delta$ c.f. \cite[Sec. 3.6]{DK}:
Specifically, the polynomials $\bar{Q}(X,0),\bar{Q}(0,\Gamma)$ are
squarefree, and $\bar{Q},\frac{\partial \bar{Q}}{\partial X},\frac{\partial \bar{Q}}{\partial \Gamma}$ have
no common solutions. Let $\bar{X}$ be the {\it toric compactification} of the
affine variety $X$ in the toric projective space $\pr_\Delta$ c.f. \cite[Sec. 3.2]{DK}. This is a smooth
compact variety.
Since the outer face
of $\Delta$ is a triangle with no interior points, it follows that $\bar{X} = X_{aff} \sqcup \pr^1$
where $X_{aff} :=\Spec( \fq[X,\Gamma,Z]/(Z^2 - \bar{Q}(X,\Gamma)))$. 

One does not need to be familiar with the exact details of the construction: the point
is simply that we have compactified the zero set in affine space of the equation
$Z^2 = \bar{Q}(X,\Gamma)$ by adding a single projective line.

\begin{definition}
  Let $P(T) \in \Z[T]$, $q$ a prime power, and $\ell$ a non-negative integer. 
  We call $P(T)$ {\rm pure of weight $\ell$} with respect to $q$ if its reciprocal zeros have complex absolute value $q^{\ell /2}$. We shall just say $P$ is a weight $\ell$ {\rm Weil polynomial}, 
  when $q$ is understood.
  \end{definition}

\begin{proposition}\label{Prop-ZetaXbar}
Let $\bar{X}$ be the smooth toric compactification of the set of affine solutions of the
equation $Z^2 = \bar{Q}(X,Y)$, as defined immediately above. Then the
zeta function of $\bar{X}$ has the form
\[
 Z(\bar{X},T) = \frac{1}{(1-T)P_2(\bar{X},T)(1 - q^2 T)},\] 
where
$P_2(\bar{X},T) \in \Z[T]$ is a Weil polynomial w.r.t. $q$ of weight
$2$, and 
\begin{equation}\label{Eqn-MiddleBettiNo}
 \deg(P_2(\bar{X},T)) = l^*(2 \Delta) - 4l^*(\Delta) - 3 -
\sum_{\Delta^\prime} (l^*(\Delta^\prime) - 1).
\end{equation}
 Here the function
$l^*$ counts lattice points in the interior of a polytope or polygon,
and the sum is over the two-dimensional faces $\Delta^\prime$ of
$\Delta$. 
\end{proposition}

\begin{proof}
The claim on the weight follows from Deligne's theorem
\cite{D}. The other claims follow from the formula for Hodge-Deligne
numbers of complex toric surfaces in \cite[Sec. 5.11(c)]{DK}, the
comparison theorem between singular cohomology and $l$-adic \'{e}tale
cohomology for the modular reduction of smooth complete varieties over
number fields, and the trace formula in $\ell$-adic \'{e}tale
cohomology.
\end{proof}

We note, but do not use, that the formula for $\deg(P_2(\bar{X},T))$ is valid for arbitrary Newton polytopes $\Delta$, 
assuming that $\pr_\Delta$ is smooth and $Z^2 - Q(X,\Gamma)$ is non-degenerate
w.r.t. $\Delta$. 

To simplify the statement and proof of the next theorem, and again to keep in line with
our implementation, we make some further restrictions: Assume that the Sylvester resultant
$\tilde{r}(\Gamma)$ is squarefree modulo $p$ of degree $d = \deg_\Gamma(r)$, 
and that for each $\gamma \in \fqbar$ with
$\tilde{r}(\gamma) = 0 \bmod{p}$, the ``missing fibre at $\Gamma = \gamma$'' in the pencil $X \rightarrow S$ has a unique double point. 

\begin{proposition}\label{Thm-FullZeta}
Definitions and assumptions as in Sections \ref{Sec-ExplicitPencil}, \ref{Sec-LMA}, \ref{Sec-P1} and the present
section. Let $\bar{X} = X \sqcup C$ where $C$ is a union of curves. We recall that
$\bar{X}$ is a compactification of the smooth surface defined by the equation
$Z^2 = \bar{Q}(X,\Gamma)$ over $\fq$, with $2g + 1:=\deg_X(\bar{Q})$ and
$h :=\deg_\Gamma(\bar{Q})$. Then $Z(C,T)$ may be
computed deterministically in $\tilde{\Oh}(g^5 h p\log(q)^3)$ bit operations. Moreover, given
$Z(C,T)$ the zeta function $Z(\bar{X},T)$ may be recovered from the first
$d - 2g$ coefficients in $P_2(X,T)$ each to $p$-adic precision modulo
$p^N$ where 
\[ N := \left\lceil \max_{0 \leq i \leq d-2g} \left\{\log_p\left(2q^i{d-2g \choose i} \right)\right\} \right\rceil.\] 
\end{proposition}

\begin{proof} Define $\bar{r}(\Gamma) := r(\Gamma) \bmod{p}$.
Let $\bar{r}(\Gamma) = \prod_{i = 1}^s
\bar{r}_i(\Gamma)$ be the irreducible factorisation and define $d_i :=\deg(\bar{r}_i)$. 
For $i = 1,\dots,s$ denote
$\gamma_i := \Gamma \in K_i := \fq[\Gamma]/\bar{r}_i(\Gamma)$. Our assumption that
each singular fibre has a unique double point implies that
$\bar{Q}(X,\gamma_i) = 
(X - \alpha_i)^2 H_i(X)$ where $H_i(\alpha_i) \ne 0$. Define
$\delta_i = -1$ if $H_i(\alpha_i)$ is a square in $K_i$, and $\delta_i := + 1$  otherwise.
Since $\bar{X} = X_{aff} \sqcup \pr^1$, it follows that 
\begin{equation}\label{Eqn-ZetaC}
 Z(C,T) =  \frac{1}{(1-T)(1-qT)}\prod_{i = 1}^s \frac{P_i(T^{d_i}) (1 + \delta_i T^{d_i})}{(1 - q^{d_i} T^{d_i})}
 \end {equation}
where $P_i(T)$ is the numerator of the zeta function of the genus $g-1$ curve $Z^2 = H_i(X)$.
Note that $Z(\bar{X},T) = Z(X,T)Z(C,T)$. From (\ref{Eqn-ZetaX}), Propositions \ref{Prop-P1} and
\ref{Prop-ZetaXbar}, and (\ref{Eqn-ZetaC}), and by noting the weights of the different
factors, we deduce the following:
\[
\begin{array}{rcl}
P_2(X,T) & = & w_2(P_2(X,T)) \prod_{i = 1}^s P_i(T^{d_i}) (1 + \delta_i T^{d_i})\\
P_2(\bar{X},T) & = & (1-qT)w_2(P_2(X,T)).
\end{array}
\]
Here $w_2(P_2(X,T))$ is the ``interesting'' weight two factor in $P_2(X,T)$. It has
degree $2g(d-1)-(d + d(2g-2)) = d - 2g$. 
The theorem now follows, using Kedlaya's algorithm \cite{KK1} to compute $Z(C,T)$ and
noting $d = \Oh(gh)$.
\end{proof}

Note that $w_2(P_2(X,T))$ satisfies the same functional equation as
$P_2(\bar{X},T)$. The sign in this functional equation is $(-1)^s$ where $s$ is the multiplicity
of $(1 + qt)$ as a factor of $P_2(\bar{X},T)$. This is unknown. However, by computing
only the first $\lfloor (d-2g)/2 \rfloor + 1$ coefficients in $P_2(X,T)$ to $p$-adic precision
modulo $p^N$ where
\begin{equation}\label{Eqn-FinalAcc}
N := \left\lceil \log_p\left( 2q^{e} {d - 2g \choose e}\right) \right\rceil,\, e:= \left\lfloor \frac{d-2g}{2} \right\rfloor
\end{equation}
one can find two possible
candidates for $P_2(\bar{X},T)$. One hopes that only one is a weight $2$ Weil polynomial!

\subsection{Computation of the full zeta function}\label{Sec-Complexity}

In this section we retain the definitions and assumptions in Sections \ref{Sec-ExplicitPencil},
\ref{Sec-LMA}, and Section \ref{Sec-Compact}. Theorem \ref{Thm-TimeSpace} and
Proposition \ref{Thm-FullZeta} together yield an algorithm for computing the full zeta function of the
compact surface $\bar{X}$, provided we can estimate the loss of precision between
the computation of the absolute Frobenius matrix $F:H_{rig}^2(X) \rightarrow H_{rig}^2(X)$
and the calculation of coefficients in the polynomial $P_2(X,T) = 
\det(1 - T q^2 F^{-\log_p(q)}| H^2_{rig}(X))$. Note that in practice one actually computes
coefficients in the polynomial $\det(T - F^{\log_p(q)}| H^2_{rig}(X))$.

\begin{theorem}\label{Thm-Varyh}
Fix a positive constant $C$ and positive integer $g$.
Assume that $\deg_X(\bar{Q}(X,\Gamma)) = 2g+1$ and that
$h := \deg_\Gamma(\bar{Q}(X,\Gamma))$ satisfies $h/p \leq C$.
Then one may compute the zeta function $Z(\bar{X},T)$ of the
compactification $\bar{X}$ of the affine surface defined by $Z^2 = \bar{Q}(X,\Gamma)$ in
$\tilde{\Oh}(h^4 p \log(q)^3)$ bit operations using $\tilde{\Oh}(h^3 p \log(q)^3)$ bits of space.
\end{theorem}

Note that the hidden constants in the Soft-Oh notation depend upon both the
genus $g$ and constant $C$.

\begin{proof}
Since $g$ is fixed and $h/p$ is bounded, the numbers $\alpha$ and $\beta$ in
Theorem \ref{Thm-Growth}, $\alpha^\prime$ in Theorem \ref{Thm-DiffSysLoss},
and consequently  $\alpha^{\prime \prime}$ and
$\beta^{\prime \prime}$ in Section \ref{Sec-NumApp} are bounded absolutely, independent
of $\bar{X}$. It follows easily from Theorem
\ref{Thm-Growth} and Propositions \ref{Prop-Eff1} and \ref{Prop-Eff2},
that the Frobenius matrix $F$ has valuation bounded below by
some absolute constant $-c$, with $c \geq 0$. We require the final $p$-adic precision to be modulo 
$p^N$ with $N$ as in the statement of Proposition \ref{Thm-FullZeta}. Notice $N = \Oh(\log(q) gh)$.
A naive analysis of the loss of accuracy during the computation of the characteristic
polynomial from the absolute Frobenius matrix shows that it certainly suffices to take $N_O = N + c \log_p(q)(d - 2g) + \ord_p((d - 2g)!)$ 
in Theorem \ref{Thm-TimeSpace}; recall $d = \Oh(gh)$. Note that condition (\ref{Eqn-GrowthNO}) is
trivially satisfied in this case since the LHS is bounded absolutely. 
The complexity estimate follows by putting this value for $N_O$ in Theorem \ref{Thm-TimeSpace}, and
noting that the resulting time/space estimate also suffices for computing the characteristic polynomial.
\end{proof}

We note that assuming Conjecture \ref{Conj-FinitePoles} is true, and that the pole order is bounded
in some manner depending only on $g$, them Theorem \ref{Thm-Varyh} holds without
the restrictions on the relative growth of $h$ and $p$. The point is that in this case one
can take $\alpha^{\prime \prime}$ and $\beta^{\prime \prime}$ to depend only on
$g$, by using the original Christol-Dwork theorem, see Note \ref{Note-Constants}.

The author has been unable to prove that the valuation of the Frobenius matrix
$F$ is bounded below by some {\it absolute} constant, i.e., bounded independent
of $g$ and $h$. If one could show this then putting $N_O = \Oh(gh \log(q))$ in Theorem
\ref{Thm-TimeSpace}, and assuming $h/p$ remains bounded as $h$ and $p$ vary,
we get the estimate $\tilde{\Oh}(g^7 h^4 p \log(q)^3)$ bit operations/$\tilde{\Oh}(g^5 h^3 p \log(q)^3)$
bits of space, for the computation of the zeta function $Z(\bar{X},T)$. Note that since the middle
betti number $d - 2g + 1$ in this case is approximately $gh$, this compares favourably
with the ``deformation algorithm''; see the end of Section \ref{Sec-Def}.

\section{Surfaces: implementation and experiments}\label{Sec-CompData}

In this section we report on a Magma (v.2.11-2) implementation of our algorithm
for the surfaces described in Sections \ref{Sec-AffineSurface} and \ref{Sec-CompactSurface} in the case of a prime field.
The experiments detailed were peformed using a dual processor Intel
Pentium 4 (3GHz with 1024Mb cache and 2Gb RAM per processor). Time and space
requirements stated are as returned by the in-built Magma function.

\subsection{Examples}

All of the examples satisfied the hypothesis in the statement of
Proposition \ref{Thm-FullZeta}, and
the local monodromy eigenvalues were as observed in Note \ref{Note-LME}.

\begin{example} \label{Ex-p17}
Let
\[ Q(X,\Gamma) := X^3 + (4\Gamma^4 + 5\Gamma^3)X 
+ \Gamma^{13} + 6\Gamma^{12} + 5\Gamma^{10} + 8\Gamma^9 + 8\Gamma^8 
+ 5\Gamma^5 +  \Gamma^4 + 5\Gamma^3 + \Gamma^2 + 1,\]
and $p:= 17$. Then the Sylvester resultant $r(\Gamma) := \Res(X,Q,\frac{\partial Q}{\partial X})$
is squarefree modulo $p$ and equals, up to a constant, the denominator $r(\Gamma)$ of
the connection matrix $b(\Gamma)/r(\Gamma)$. Both polynomials have degree $d := 2 \times
13 = 26$. The genus of the generic fibre is $g := 1$. The space $H^2_{rig}(X)$ associated to the open surface $Z^2 = Q(X,\Gamma), r(\Gamma) \ne 0 \bmod{p}$ has dimension $(d-1)2g = 50$. The space
$H^2_{rig}(\bar{X})$ associated with the smooth toric compactification has dimension $d - 2g + 1 = 
25$. Computing a matrix for the Frobenius map $F:H^2_{rig}(X) \rightarrow H^2_{rig}(X)$ to precision
modulo $p^{18}$, we recovered two possible choices for the polynomial $\det(T - F|H^2_{rig}(\bar{X}))$.
Only one was the reciprocal of a weight $2$ Weil polynomial w.r.t. $17$. Specifically, 
$P_2(\bar{X},T) = (1 - 17T)^2 (1 + 17T) R(T)$ where $R(T)$ is the irreducible polynomial
{\small
\[
1+2^{3}T^{1}+2^{1}3^{2}17^{1}T^{2}+2^{4}17^{1}19^{1}T^{3}+2^{2}3^{2}5^{1}17^{2}T^{4}
+3^{1}17^{3}23^{1}T^{5}+17^{4}23^{1}T^{6}\]
\[-2^{1}5^{1}17^{5}T^{7}+3^{4}17^{6}T^{8}+2^{2}5^{1}7^{1}17^{7}T^{9}+2^{1}17^{8}191^{1}T^{10}
+2^{4}13^{1}17^{9}T^{11}+2^{1}17^{10}191^{1}T^{12}\]
\[+2^{2}5^{1}7^{1}17^{11}T^{13}+
3^{4}17^{12}T^{14}-2^{1}5^{1}17^{13}T^{15}+17^{14}23^{1}T^{16}+3^{1}17^{15}23^{1}T^{17}+2^{2}3^{2}5^{1}17^{16}T^{18}\]
\[+2^{4}17^{17}19^{1}T^{19}+2^{1}3^{2}17^{19}T^{20}+2^{3}17^{20}T^{21}+17^{22}T^{22}.
\]
}
The Hodge numbers defined a polygon called the Hodge polygon which lies below the Newton
polygon of $P_2(X,T)$. In this case, the Hodge numbers are $2,21,2$ which explains the
high divisilbility of the coefficients by powers of $p$, c.f. \cite[Remark 1.6.4]{AKR}.

The computation is provably correct under no additional hypothesis. It took just under $23$ hours and
$13$ minutes, and required just under 1.312 Gbytes of memory. We note that over half
the time required was taken computing $r(\Gamma)$-expansions of the elements
in the our relative Frobenius matrix. This was necessary to ensure the input to the second
stage of the algorithm was in the appropriate form.  

Under Conjecture \ref{Conj-FinitePoles}. this example required just under 2 hours 36 minutes, 
and 216Mbytes of memory --- the pole order appears to be $39$.
\end{example}

\begin{example}\label{Ex-p5}
Let
\[ Q(X,\Gamma) := X^3 + (\Gamma^{13} + 3\Gamma^3 + 1)X + \Gamma^{31} 
+ 2\Gamma^{15} + 4\Gamma^8 + 3\Gamma^3 + 2\Gamma + 1,\]
and $p:=5$. So $d = 2 \times 31 = 62$, $g :=1$, $\dim(H^2_{rig}(X)) = 2g(d-1) = 121$,
$\dim(H^2_{rig}(\bar{X})) = d - 2g + 1 = 61$. The characteristic polynomial
$\det(T - F|H^2_{rig}(X))$ was computed modulo $p^{55}$.
We found $P_2(\bar{X},T) = 
(1-5T)^2(1 + 5T)R(T)$ where $R(T)$ is the irreducible integer polynomial
{\small
\[
1-5^{1}T^{1}+2^{1}3^{1}5^{1}T^{2}-2^{1}3^{1}5^{1}T^{3}-2^{3}5^{2}T^{4}+5^{2}89^{1}T^{5}-2^{1}5^{2}409^{1}T^{6}+2^{1}5^{3}7^{1}47^{1}T^{7}\]
\[-3^{1}5^{3}617^{1}T^{8}-5^{4}37^{2}T^{9}+5^{5}7^{1}727^{1}T^{10}-3^{1}5^{6}11^{1}97^{1}T^{11}+5^{8}11^{1}53^{1}T^{12}+5^{9}7^{1}73^{1}T^{13}\]
\[-5^{9}7^{2}43^{1}T^{14}+5^{11}11^{1}61^{1}T^{15}-2^{1}5^{11}829^{1}T^{16}-2^{1}5^{12}677^{1}T^{17}+2^{1}3^{1}5^{14}79^{1}T^{18}\]
\[-5^{14}53^{1}89^{1}T^{19}+2^{1}5^{15}1777^{1}T^{20}-5^{18}37^{1}T^{21}+2^{1}3^{1}5^{17}11^{1}T^{22}+5^{20}137^{1}T^{23}-3^{1}5^{20}7^{2}T^{24}\]
\[+2^{3}3^{2}5^{20}29^{1}T^{25}+5^{21}367^{1}T^{26}-5^{23}7^{1}23^{1}T^{27}+2^{1}3^{2}5^{23}53^{1}T^{28}-2^{3}3^{1}5^{24}7^{1}13^{1}T^{29}\]
\[+2^{1}3^{2}5^{25}53^{1}T^{30}-5^{27}7^{1}23^{1}T^{31}+5^{27}367^{1}T^{32}+2^{3}3^{2}5^{28}29^{1}T^{33}-3^{1}5^{30}7^{2}T^{34}+5^{32}137^{1}T^{35}\]
\[+2^{1}3^{1}5^{31}11^{1}T^{36}-5^{34}37^{1}T^{37}+2^{1}5^{33}1777^{1}T^{38}-5^{34}53^{1}89^{1}T^{39}+2^{1}3^{1}5^{36}79^{1}T^{40}\]
\[-2^{1}5^{36}677^{1}T^{41}-2^{1}5^{37}829^{1}T^{42}+5^{39}11^{1}61^{1}T^{43}-5^{39}7^{2}43^{1}T^{44}+5^{41}7^{1}73^{1}T^{45}\]
\[+5^{42}11^{1}53^{1}T^{46}-3^{1}5^{42}11^{1}97^{1}T^{47}+5^{43}7^{1}727^{1}T^{48}-5^{44}37^{2}T^{49}-3^{1}5^{45}617^{1}T^{50}\]
\[+2^{1}5^{47}7^{1}47^{1}T^{51}-2^{1}5^{48}409^{1}T^{52}+5^{50}89^{1}T^{53}-2^{3}5^{52}T^{54}-2^{1}3^{1}5^{53}T^{55}+2^{1}3^{1}5^{55}T^{56}\]
\[-5^{57}T^{57}+5^{58}T^{58}.
\]
 }
The Hodge numbers in this case are $5,51,5$.

The computation is provably correct only under Conjecture \ref{Conj-FinitePoles} --- the pole
order appears to be 31. It took
13 hours and 577 seconds, and required just under 834 Mbytes of memory.
\end{example}

\begin{example}\label{Ex-p11}
 Let
\[ Q(X,\Gamma) := X^5 + 4X^3 + (4\Gamma^2 + 4\Gamma + 8)X + \Gamma^7 + 5\Gamma^6 + 1\]
and $p:=11$.  So $d = 4 \times 7 = 28$, $g :=2$, $\dim(H^2_{rig}(X)) = 2g(d-1) = 108$,
$\dim(H^2_{rig}(\bar{X})) = d - 2g + 1 = 25$. The characteristic polynomial
$\det(T - F|H^2_{rig}(X))$ was computed modulo $p^{19}$.
We found $P_2(\bar{X},T) = 
(1-11T)^3R(T)$ where $R(T)$ is the irreducible integer polynomial
{\small
\[
1+19^{1}T^{1}+17^{2}T^{2}+2^{5}11^{2}T^{3}+7^{3}11^{2}T^{4}+3^{1}11^{3}103^{1}T^{5}
+2^{1}3^{1}11^{4}41^{1}T^{6}+3^{2}11^{5}23^{1}T^{7}\]
\[+11^{6}151^{1}T^{8}+2^{5}5^{1}11^{7}T^{9}+2^{2}11^{8}47^{1}T^{10}+3^{1}11^{9}59^{1}T^{11}
+2^{2}11^{10}47^{1}T^{12}+2^{5}5^{1}11^{11}T^{13}\]
\[+11^{12}151^{1}T^{14}+3^{2}11^{13}23^{1}T^{15}+2^{1}3^{1}11^{14}41^{1}T^{16}
+3^{1}11^{15}103^{1}T^{17}+7^{3}11^{16}T^{18}\]
\[+2^{5}11^{18}T^{19}+11^{18}17^{2}T^{20}+11^{20}19^{1}T^{21}+11^{22}T^{22}.\]
}
The Hodge numbers are $2,21,2$ in this case.

The computation is provably correct only under Conjecture \ref{Conj-FinitePoles} --- the pole
order appears to be $21$. It took just under 14 hours 36 minutes and required 4.41Gbytes of memory.
\end{example}

We note that use of Method 2 (Section \ref{Sec-Method2}) rather than Method 1 (Section
\ref{Sec-Method1}) significantly reduces the space requirement; however, we do not have 
provable precision loss bounds for Method 2. If one is satisfied with {\it plausible} rather than
{\it provable} output, larger examples may be computed.

\subsection{Calculation of precisions required}\label{Sec-ExactPrecision}

We now address the delicate problem of minimizing the amount of precision
one needs to carry through the algorithm to obtain an answer which is provable
correct (possibly assuming Conjecture \ref{Conj-FinitePoles}).

Fix a positive integer $N_3$ and suppose that we wish to compute a $p^{N_3}$-approximation
to a matrix for $F$ acting on $H^2_{rig}(X)$. 
We will discuss the choice of $N_3$ later in this section.
Let
$B_{2g,p}$ be as in (\ref{Eqn-Bmp}); in particular, for $p \geq 2g$ we have
$B_{2g,p} = 2g-1$. Recall that $2g + 1 :=\deg_X(Q)$.

Define $x_{fin}$ and $N_{2,fin}$ to be the smallest integer solutions to the inequalities:
\[
\begin{array}{rcl}
\lfloor x_{fin}/p \rfloor - \lfloor \log_p(2x_{fin} + 1) \rfloor & \geq & N_{2,fin}\\
N_{2,fin} - (2 B_{2g,p} + 2g) \lfloor \log_p(x_{fin}) \rfloor & \geq & N_3.
\end{array}
\]
More precisely, let $x_{fin}$ be the smallest integer solution to
\[ 
\lfloor x_{fin}/p \rfloor - \lfloor \log_p(2x_{fin} + 1)\rfloor - (2 B_{2g,p} + 2g) \lfloor \log_p(x_{fin}) \rfloor 
 \geq  N_3,
\]
and define $N_{2,fin}$ in the obvious way.
Applying Proposition \ref{Prop-Eff1} and
Theorem \ref{Thm-Growth}, and recalling that we have nilpotent monodromy around
the roots of $r(\Gamma)$, one sees the following: It is enough to compute the coefficients
$f_k(\Gamma)$ with $k < 0$ in the $r(\Gamma)$-adic expansion of the entries of $F(\Gamma)$
for $|k| \leq x_{fin}$, and to compute these with $p$-adic precision ``modulo $p^{N_{2,fin}}$''.
The point is that for any basis form $b_{ik}(\Gamma)$ (Section \ref{Sec-LMA}),
the coefficients in ``$(i,k)$th column'' of the matrix for $F$ are given by applying the reduction 
algorithm from the proof of Theorem \ref{Thm-FGEdR} to the image form
$F(\Gamma)b_{ik}(\Gamma^p)$; but the reduced 
form to $p$-adic precision ``modulo $p^{N_3}$'' is not affected by
coefficients $f_k(\Gamma)$ for $k$ negative with $|k| > x_{fin}$. 
Define $N_{\Gamma,fin} := \deg(r)x_{fin}$.

One can argue in a similar manner to determine which coefficients $f_k(\Gamma)$
for $k \geq 0$ in the $r(\Gamma)$-adic expansions of entries in $F(\Gamma)$ must be computed, and to what precision. We return to this shortly, but let us say that we have determined suitable integers 
$x_{inf}$ and $N_{2,inf}$, and defined $N_{\Gamma,\inf} :=\deg(r)x_{inf}$.

Define $N_2 := \max\{N_{2,fin},N_{2,inf}\}$. We need to compute the coefficients $f_k(\Gamma)$ in the
$r$-adic expansion of entries in the global Frobenius matrix $F(\Gamma)$ for
$-x_{fin} \leq k < x_{inf}$ with $p$-adic precision ``modulo $p^{N_2}$''. Define
$N_\Gamma := N_{\Gamma,fin} + N_{\Gamma,inf}$. Since there is no loss of accuracy
during the analytic continuation stage (Section \ref{Sec-GlobalDef}), it is enough to compute
a $p^{N_2}$-approximation to  the local Frobenius matrix $F(\Gamma)$ modulo $\Gamma^{N_{\Gamma}}$. Using the method in Section
\ref{Sec-Method1}, equation (\ref{Eqn-M1loss}) tells us we must perform
the local calculation itself to $p$-adic precision  ``modulo $p^{N_1}$'', where
\[ N_1 := N_2 + (3B_{2g,p} + 1) \lfloor \log_p(N_\Gamma)\rfloor - B_{2g,p} + \min\{\ord_p(F(0)),0\}.\]
Note that $\ord_p(F(0)) \geq 0$ when $p \geq 2g$.
Our algorithm begins by computing a $p^{N_1}$-approximation to the matrix $F(0)$; 
see \cite{KK1} for an analysis of the loss of accuracy during this
initial computation.

We return to the question of determining $N_{2,inf}$, $x_{\inf}$ and $N_{\Gamma,inf}$. One can do this
via an analogous system of inequalities to those above, using Proposition
\ref{Prop-Eff2} and Theorem \ref{Thm-Growth} combined with the general
estimates for $\alpha$ and $\beta$ derived from Note \ref{Note-Constants}. 
The problem is that since the local monodromy around the point at infinity is not nilpotent, the
constants $\alpha$ and $\beta$ are rather large. The get around this, the author wrote
a short computer program which calculated more careful bounds on the growth of the 
coefficients in the uniform part of the local solution matrix around the point at infinity. 
In the notation of Lemma \ref{Lem-DGS}, the author
computed a lower convex function $a_1(i)$ such that $\ord_p(Y_i) \geq - a_1(i)$ for all
$i \geq 1$. The function $a_1(i)$ depended explicitly on the local monodromy eigenvalues
at infinity and $p$; the time required to compute $a_1(i)$ grew as $\log_p(i)$ with $i \geq 1$.
Here are brief details: For eigenvalues in the interval $[0,1)$ use \cite[Lines 8-9,19, Page 196]{DGS};
for general prepared eigenvalues, use the proof of Lemma \ref{Lem-DGS}, but
compute a tighter lower bound on ``$\ord_p(\tilde{\h})$'' via the proof of Lemma \ref{Lem-Naive}, and
use the inequality ``$\ord_p(Y_i) \geq \ord_p(\tilde{Y}_{i + 2\Delta}) + \ord_p(\tilde{\h})$''. 
The function $a_1(i)$ was fed as input to the analysis in the proof of Theorem \ref{Thm-Growth}, to
yield a better function $a(\ell)$, say, which could be used on the righthand-side on the
statement of the theorem. With this more refined function, one takes $x_{inf}$ and $N_{2,inf}$ to be the smallest integer solutions to the inequalities:
\[
\begin{array}{rcl}
(y_{inf} + 1) - \lfloor \log_p(p(2y_{inf} + 1)) \rfloor & \geq & N_{2,inf}, \left(y_{inf} := 
\left\lfloor \frac{2 x_{inf} -  (p-1)(\delta  -1)}{2p(\delta + \delta^\prime - 1)} \right\rfloor \right)\\
N_{2,inf} - a(x_{inf} \deg(r)) & \geq & N_3.
\end{array}
\]
See (\ref{Eqn-Delta}) for the definitions of  the numbers $\delta$ and $\delta^\prime$.

When the author assumed Conjecture \ref{Conj-FinitePoles} was true, he did not perform
the calculation in the preceding paragraph, but instead defined $N_{2,inf} := N_{2,fin}$ and $N_{\Gamma,inf} := 100$.

We require that the characteristic polynomial $\det(T - F| H^2_{rig}(X))$ be computed
modulo $p^N$, with $N$ as in equation (\ref{Eqn-FinalAcc}). If $\ord_p(F) \geq 0$, then one
can take $N_3 := N$; this was the case in Example \ref{Ex-p17}. If $\ord_p(F) < 0$, there
may be some loss of accuracy during the computation of the characteristic polynomial. The author
had an {\it ad hoc} solution to this problem: Specifically, it was observed in practice that even
when $\ord_p(F) < 0$, some small power of $F$ had non-negative or even positive valuation. By examining the valuation of powers of $F$, and using the formula $P_2^\prime(X,T)/P_2(X,T) = 
-\sum_{k = 1}^\infty \Tr(F^k) T^{k-1}$, one can deduce explicit bounds on the loss of precision. This enabled
the author to establish usable and provable precision loss bounds during the calculation of the
characteristic polynomial; however, when the initial computation revealed $\ord_p(F) < 0$, one
did need to rerun the computation with an increased value for $N_3$ to get a provably correct
answer. 

The parameters $[N,N_3,N_{2,fin},N_{2,inf},N_1;N_{\Gamma,fin},N_{\Gamma,inf}]$ in the examples
were set as follows: in Example \ref{Ex-p17}, $[18,18,26,56,67;12376,16692]$ unconditionally and
$[18,18,26,26,37;12376,100]$ under Conjecture \ref{Conj-FinitePoles}; in Example
\ref{Ex-p5}, $[55,60,72,72,95;23560,100]$; in Example \ref{Ex-p11}, $[19,25,45,45,72;14476,100]$.


\begin{thebibliography}{99}

\bibitem{AKR} T.G. Abbot, K. Kedlaya, D. Roe, Bounding Picard numbers of surfaces using
$p$-adic cohomology, to appear in Proc. Arith. Geom. and Coding Theory (Luminy 2005).
See www-math.mit.edu/$\sim$kedlaya/papers/ and arxiv.org/abs/math.NT/0601508

\bibitem{AA} A. Adolphson, An index theorem for $p$-adic differential operators, Trans.
A.M.S. Vol. 216, (1976), 279-293.


\bibitem{BC} F. Baldassarri and B. Chiarellotto, Algebraic versus rigid cohomology with
logarithmic coefficients, 11-150, Barsotti Symposium on Algebraic Geometry,V. Cristante,
W. Messing (Eds), Academic Pres, 1994.

\bibitem{PB} P. Berthelot, Cohomologie rigide et cohomologie rigide \`{a} supports propres,
Premi\`{e}re partie (version provisoire 1991), Pr\'{e}publication 96-03, Institut de Recherche
Math\'{e}matique de Rennes, 1996. 

\bibitem{CD} G. Christol and B. Dwork, Effective $p$-adic bounds at regular
singular points, Duke. Math. 62, (1991), 689-720.

\bibitem{DNC} D.N. Clark, A note on the $p$-adic  convergence of solutions of linear
differential equations, Proc. A.M.S. Vol. 17, No.1, 262-269.

\bibitem{DK} V.I. Danilov and A.G. Khovanski\v{i}, Newton polyhedra and an algorithm for computing
Hodge-Deligne numbers, Math. USSR Izvestiya Vol. 29 (1987), No.2.

\bibitem{CLO} D. Cox, J. Little and D. O'Shea, Ideals, Varieties, and Algorithms,
2nd Edition, Springer UTM, 1997.

\bibitem{D} P. Deligne, La conjecture de Weil: I, Pub. I.H.E.S. 43, (1974), 273-307.

\bibitem{DGS} B. Dwork, G, Gerotto, F.J. Sullivan, An Introduction to G-functions, Annals
of Mathematical Studies 133, Princeton University Press, 1994.

\bibitem{vGG} J. von zur Gathen and J. Gerhard, Modern Computer Algebra, Cambridge University
Press, 1999.

\bibitem{G2} R. Gerkmann, Relative rigid cohomology and point counting on families of
elliptic curves, preprint 2005, available with Magma implementation at: www.mathematik.uni-mainz.de/$\sim$gerkmann/.

\bibitem{HH} H. Hubrechts, Point counting on families of hyperelliptic curves,
preprint 2005. See http://arxiv.org/math.NT/0601438.

\bibitem{NK} N.M. Katz, Nilpotent connections and the monodromy theorem: applications
of a result of Turrittin, Publ. IHES 39 (1970), 175-232.

\bibitem{NKDwork} N. Katz, Travaux de Dwork, S\'{e}minar Bourbaki 24e ann\'{e}e, 1971/72, No. 409, 
167-200.

\bibitem{NK2} N.M Katz, An overview of Deligne's proof of the Riemann hypothesis
for varieties over finite fields, A.M.S. Proc. Symp. Pure Math, Vol. 28, (1976), 275-305.

\bibitem{KO} N.M. Katz and T. Oda, On the differentiation of De Rham cohomology
classes with respect to parameters, J. Math. Kyoto Univ. 8-2, (1968), 199-213.

\bibitem{KK1} K. Kedlaya, Counting points on hyperelliptic curves using Monsky-Washnitzer
cohomology, J. Ramanujan Math. Soc., Vol. 16, (2001), 323-338. 

\bibitem{KKsur} K. Kedlaya, Computing zeta functions via $p$-adic cohomology,
ANTS 2004, D.A. Buell (ed), LNCS 3076, (2004), 1-17.

\bibitem{KKII} K. Kedlaya, Fourier transforms and $p$-adic ``Weil II'', preprint. See
http://arxiv.org/math.NT/0210149.

\bibitem{LSH} A.G.B. Lauder, Counting solutions to equations in many variables over
finite fields, Foundations of Computational Mathematics Vol. 4. No. 3, (2004), 221-267.

\bibitem{LFFA} A.G.B. Lauder, Rationality and meromorphy of zeta functions, Finite Fields and
Their Applications, Vol. 11, (2005), 491-510.

\bibitem{LJNTB} A.G.B Lauder, Rigid cohomology and $p$-adic point counting, J. Th. Nom. Bordeaux
Vol 17 No. 1, (2005), 169-180.

\bibitem{LW} A.G.B. Lauder and D. Wan, Counting points on varieties over finite fields of
small characteristic, Algorithmic number theory: lattices, number fields, curves and
cryptography (ed. J.P. Buhler, P. Stevenhagen), MSRI Pub. 44, to appear.

\bibitem{JP} J. Pila, Frobenius maps of abelian varieties and finding roots of unity
in cyclotomic fields, Math. Comp. 55, (1990), 745-763.

\bibitem{vPS} M. van der put and M. Singer, Galois Theory of Linear Differential Equations,
Series of comprehensive studies in mathematics Vol. 328,
Springer, 2003.

\bibitem{RS} R. Schoof, Elliptic curves over finite fields and the computation of square
roots modulo $p$, Math. Comp. 44, (1985), 483-494.

\bibitem{MS} M. Setoyanagi, Note on Clark's theorem for $p$-adic convergence,
Proc. A.M.S. Vol. 125, No.3, 717-721.

\bibitem{NT} N. Tsuzuki, Bessel $F$-isocrystals and an algorithm for computing Kloosterman sums, 
preprint 2003.

\bibitem{NTss} N. Tsuzuki, On base change theorem and coherence in rigid cohomology, Documenta
Math., Extra Volume Kato, (2003), 891-918.

\bibitem{DW} D. Wan, Algorithmic theory of zeta functions, MSRI Pub. 44, to appear.



\end{thebibliography}
\end{document}